\documentclass[final]{siamart220329}

\usepackage{showkeys}

\usepackage[show]{ed}

\usepackage{tikz}
\usepackage{pgf}
\usetikzlibrary{arrows}
\usetikzlibrary{calc}
\usetikzlibrary{decorations.pathmorphing}
\usepackage{xfrac}    

\usepackage{faktor}
\usepackage{soul}
\setstcolor{red}
\setulcolor{red}

\usepackage{rotating} 

\usepackage{mathrsfs}
\DeclareMathAlphabet{\mathpzc}{OT1}{pzc}{m}{it}

\usepackage{color}
\usepackage{amsmath}
\usepackage{graphicx}
\usepackage{graphics}
\usepackage{float}
\usepackage{amsfonts}
\usepackage{amssymb}%
\usepackage{latexsym}
\usepackage{psfrag}
\usepackage{subfigure}
\usepackage{accents}

\newtheorem{assumption}[theorem]{Assumption}
\newtheorem{remark}[theorem]{Remark}
\newtheorem{example}[theorem]{Example}

\newcommand{\cF}{{\cal F}}
\newcommand{\cH}{{\cal H}}
\newcommand{\cI}{{\cal I}}

\newcommand{\cP}{{\cal P}}

\newcommand{\cR}{{\cal R}}
\newcommand{\cS}{{\cal S}}
\newcommand{\cT}{{\cal T}}

\newcommand{\cV}{{\cal V}}

\newcommand{\cN}{{\cal N}}

\newcommand{\bx}{\mathbf{x}}
\newcommand{\br}{\mathbf{r}}

\newcommand{\bu}{\mathbf{u}}

\newcommand{\bff}{\mathbf{f}}

    \newcommand\quotient[2]{
        \mathchoice
            {
                \text{\raise1ex\hbox{$#1$}\Big/\lower1ex\hbox{$#2$}}%
            }
            {
                #1\,/\,#2
            }
            {
                #1\,/\,#2
            }
            {
                #1\,/\,#2
            }
    }

\newcommand{\ri}{{\rm i}}
\newcommand{\rd}{{\rm d}}

\newcommand{\beq}{\begin{equation}}
\newcommand{\eeq}{\end{equation}}
\newcommand{\beqs}{\begin{equation*}}
\newcommand{\eeqs}{\end{equation*}}
\newcommand{\bit}{\begin{itemize}}
\newcommand{\eit}{\end{itemize}}
\newcommand{\ben}{\begin{enumerate}}
\newcommand{\een}{\end{enumerate}}
\newcommand{\bal}{\begin{align}}
\newcommand{\eal}{\end{align}}
\newcommand{\bals}{\begin{align*}}
\newcommand{\eals}{\end{align*}}
\newcommand{\bse}{\begin{subequations}}
\newcommand{\ese}{\end{subequations}}
\newcommand{\bpr}{\begin{proposition}}
\newcommand{\epr}{\end{proposition}}
\newcommand{\bre}{\begin{remark}}
\newcommand{\ere}{\end{remark}}
\newcommand{\bpf}{\begin{proof}}
\newcommand{\epf}{\end{proof}}
\newcommand{\ble}{\begin{lemma}}
\newcommand{\ele}{\end{lemma}}
\newcommand{\bco}{\begin{corollary}}
\newcommand{\eco}{\end{corollary}}
\newcommand{\bex}{\begin{example}}
\newcommand{\eex}{\end{example}}
\newcommand{\bth}{\begin{theorem}}
\newcommand{\enth}{\end{theorem}}

\newcommand{\Rea}{\mathbb{R}}
\newcommand{\Com}{\mathbb{C}}

\def\XXint#1#2#3{{\setbox0=\hbox{$#1{#2#3}{\int}$}
     \vcenter{\hbox{$#2#3$}}\kern-.5\wd0}}

\usepackage{hyperref}
\counterwithin{figure}{section}
\definecolor{myblue}{rgb}{0,0,0.6}
\newcommand*{\N}[1]{\left\|#1\right\|}

\allowdisplaybreaks[4]

\newcommand{\tfa}{\text{ for all }}
\newcommand{\tfor}{\text{ for }}

\newcommand{\tor}{\text{ or }}
\newcommand{\tif}{\text{ if }}

\newcommand{\ton}{\text{ on }}
\newcommand{\tas}{\text{ as }}
\newcommand{\tand}{\text{ and }}
\newcommand{\tst}{\text{ such that }}
\newcommand{\tfind}{\text{ find }}

\newcommand{\vertiii}[1]{{\left\vert\kern-0.25ex\left\vert\kern-0.25ex\left\vert #1
    \right\vert\kern-0.25ex\right\vert\kern-0.25ex\right\vert}}

\definecolor{jwcol}{RGB}{27, 137, 18}  

\definecolor{dalcol}{rgb}{0.8,0,0}

\definecolor{escol}{rgb}{0,0,0.8}
\definecolor{estcol}{rgb}{0,0.5,0}
\definecolor{esnewcol}{rgb}{0,0.5,0}

\newcommand{\supp}{{\rm supp}}

\newcommand{\Ccom}{C_{\rm com}}
\newcommand{\sigmao}{\sigma_{H^1}}
\newcommand{\sigmat}{\sigma_{L^2}}

\newcommand{\Csol}{C_{\rm sol}}

\newcommand{\Ccont}{{C_{\rm cont}}}

\newcommand{\bV}{\mathbf{V}}

\newcommand{\bW}{\mathbf{W}}

\newcommand{\matrixI}{{\mathsf{I}}}

\newcommand{\fine}{\cV_h}
\newcommand{\coarse}{\cV_0}

\newcommand{\subdomain}{\cV_{\ell,h}}

\newcommand{\matrixD}{\mathsf{D}}
\newcommand{\matrixA}{\mathsf{A}}

\newcommand{\matrixAzero}{\matrixA_{0}}

\newcommand{\matrixS}{\mathsf{S}}
\newcommand{\matrixR}{\mathsf{ R}}
\newcommand{\matrixM}{\mathsf{ M}}
\newcommand{\matrixB}{\mathsf{ B}}

\newcommand{\Ccoarse}{C_{\rm coarse}}

\newcommand{\cJ}{{\mathcal J}}

\newcommand{\chiellg}{\chi_\ell^{>}}
\newcommand{\chiell}{\chi_\ell}
\newcommand{\Hilbert}{\cV}

\newcommand{\coeffA}{A}
\newcommand{\coeffB}{B}
\newcommand{\coeffc}{c}

\newcommand{\MAP}{\cF}
\newcommand{\LF}{\MAP}
\newcommand{\CT}{\cT}
\newcommand{\CP}{\cP}
\newcommand{\CI}{\cI}
\newcommand{\hK}{\widehat{K}}
\newcommand{\element}{K}

\newcommand{\Hot}{H^1_{\widetilde{0}}}
\newcommand{\Hcoarse}{H_{\rm coarse}}

\newcommand{\Ccurve}{C_{\rm curved}}
\newcommand{\Csuper}{C_{\rm super}}
\newcommand{\Cint}{C_{\rm int}}

\makeatletter
\newcommand{\settheoremtag}[1]{
  \let\oldthetheorem\thetheorem
  \renewcommand{\thetheorem}{#1}
  \g@addto@macro\endtheorem{
    \addtocounter{theorem}{-1}
    \global\let\thetheorem\oldthetheorem}
  }
\makeatother

\usetikzlibrary{positioning}
\usepackage{lscape}

\definecolor{jeffColor}{RGB}{102, 0, 204}

\title{
Convergence theory for two-level hybrid Schwarz preconditioners for high-frequency Helmholtz problems
}

\author{
J.~Galkowski\thanks{Department of Mathematics, University College London, 25 Gordon Street, London, WC1H 0AY, UK,   \tt J.Galkowski@ucl.ac.uk}
\and
E.~A.~Spence\thanks{Department of Mathematical Sciences, University of Bath, Bath, BA2 7AY, UK, \tt E.A.Spence@bath.ac.uk }
}

\date{\today}

\begin{document}
\maketitle

\begin{abstract}
We give a novel convergence theory for two-level hybrid Schwarz domain-decomposition (DD) methods for finite-element discretisations of the high-frequency Helmholtz equation.
This theory gives sufficient conditions for the preconditioned matrix to be close to the identity, and 
covers DD subdomains of arbitrary size, arbitrary absorbing layers/boundary conditions on both the global and local Helmholtz problems,
and coarse spaces not necessarily related to the subdomains.

The assumptions on the coarse space are satisfied by the approximation spaces using problem-adapted basis functions that have been recently analysed 
as coarse spaces for the Helmholtz equation, 
as well as all spaces in which the Galerkin solutions are known to be quasi-optimal  via a Schatz-type argument.

As an example, we apply this theory when the coarse space consists of piecewise polynomials; these are then the first rigorous convergence results about a two-level Schwarz preconditioner applied to the high-frequency Helmholtz equation with a coarse space that does not consist of problem-adapted basis functions.
\end{abstract}

\section{Introduction}

\subsection{Context and motivation}

Coarse grids are 
the key to parallel scalability of domain-decomposition (DD) methods for self-adjoint positive-definite problems (such as Laplace's equation); see, e.g., \cite{ToWi:05, Wi:09}, \cite[Chapter 4]{DoJoNa:15}. 
However, the design of practical coarse spaces (with associated theory) for high-frequency wave problems, such as the high-frequency Helmholtz equation, is a longstanding open problem 
(see, e.g., the recent computational study \cite{BoDoJoTo:21} and the references therein). 

For the Helmholtz equation with fixed wavenumber $k$, \cite{CaWi:92} analysed two-level additive Schwarz methods, 
and \cite{GSV1} performed the corresponding analysis when $k$ is complex valued with large $|k|$ and sufficiently large imaginary part.

For the Helmholtz equation with $k$ real and large, the four recent papers \cite{HuLi:24, LuXuZhZo:24, MaAlSc:24, FuGoLiWa:24} all analyse hybrid Schwarz methods (where the coarse and subdomain solves are combined in a multiplicative way) with coarse spaces consisting of problem-adapted basis functions (coming from solving Helmholtz problems on subsets of the domain).

The present paper gives a novel convergence theory for two-level hybrid Schwarz methods under very general assumptions (e.g., the subdomains can have arbitrary size and boundary conditions). This convergence theory has three immediate applications.
\begin{enumerate}
\item 
The coarse spaces in  \cite{HuLi:24, LuXuZhZo:24, MaAlSc:24, FuGoLiWa:24} all satisfy the assumptions in the theory of the present paper; 
 therefore our theory gives results about hybrid Schwarz preconditioners using these coarse spaces, complementary to those in  \cite{HuLi:24, LuXuZhZo:24, MaAlSc:24, FuGoLiWa:24}. Indeed, although the method in the present paper and the methods in \cite{HuLi:24, LuXuZhZo:24, MaAlSc:24, FuGoLiWa:24} are all hybrid, the particular hybrid method considered here is slightly different. 
However, the theory here covers arbitrary absorbing layers/boundary conditions on both the global and local Helmholtz problems, while the analyses in \cite{HuLi:24, LuXuZhZo:24, MaAlSc:24, FuGoLiWa:24} are all for global impedance boundary conditions and local Dirichlet (with no absorbing layers) or local impedance boundary conditions.
Furthermore, while \cite{HuLi:24, LuXuZhZo:24} prove results about the field of values of the preconditioned matrix, the theory here proves the stronger result that the preconditioned matrix is close to the identity.
\item The theory in the present paper applies to problem-adapted approximation spaces used to solve the Helmholtz equation that have 
not yet been used as coarse spaces, e.g., the multiscale space of \cite{ChHoWa:23}.
\item We use the theory in the present paper to prove the first convergence results about piecewise-polynomial coarse spaces for the high-frequency Helmholtz equation.
We highlight that the recent computational study \cite{BoDoJoTo:21} comparing coarse spaces found piecewise-polynomial coarse spaces to be competitive -- from the point of view of number of GMRES iterations -- with coarse spaces involving problem-adapted basis functions. Piecewise-polynomial coarse spaces also have the advantage that they are much cheaper to use than those involving problem-adapted basis functions. 
\een

\subsection{Informal statement of the main result}

\subsubsection{The global and local Helmholtz problems}

Let $\Omega\subset\Rea^d$ be bounded Lipschitz domain, and let $\Hilbert$ equal $H^1(\Omega)$, possibly with zero Dirichlet boundary conditions on part of $\partial\Omega$. 
Let
\beq\label{eq:sesqui_intro}
a(u,v):= \int_\Omega k^{-2} (\coeffA \nabla u)\cdot\overline{\nabla v} + k^{-1} (\coeffB\cdot\nabla u)\overline{v} - c^{-2} u\overline{v} - \ri k^{-1} \int_{\partial\Omega} \theta u\overline{v}
\eeq
be a sesquilinear form on $\Hilbert$, with $\coeffA\in L^\infty(\Omega;\Com^{d\times d})$, $\coeffB\in L^\infty(\Omega;\Com^{d})$, and $\coeffc \in L^\infty(\Omega;\Com)$. 
We are primarily interested in solving Helmholtz problems that approximate scattering by either a penetrable obstacle (modelled by variable coefficients $\coeffA$, $\coeffB$, and $\coeffc$ in the interior of $\Omega$) and/or an impenetrable obstacle (via $\partial\Omega$ not being connected).
Truncation of the unbounded domain exterior to the scatterer by a perfectly-matched layer (PML) \cite{Be:94, CoMo:98}, other absorbing layers (such as a complex-absorbing potential 
\cite{RiMe:93}), or an impedance boundary condition can all be written in the form of \eqref{eq:sesqui_intro} via appropriate choices of $\coeffA,\coeffB,\coeffc$, and $\theta$. 

Let $\{\Omega_\ell\}_{\ell=1}^N$ be an overlapping cover of $\Omega$, and let $\{\chiell\}_{\ell=1}^N$ be a partition of unity subordinate to  $\{\Omega_\ell\}_{\ell=1}^N$. Let $\{\chiellg\}_{\ell=1}^N$ be a second set of cut-off functions with $\chiellg$ supported in $\Omega_\ell$ and $\chiellg\equiv 1$ on $\supp\chiell$; i.e., $\chiellg$ is ``bigger than" $\chiell$, hence the notation,  and $\{\chiellg\}_{\ell=1}^N$  is \emph{not} a partition of unity.
We highlight that the existence of such $\chiellg$ means that the overlap of the subdomains must be sufficiently large; the role of $\chiellg$ is discussed in 
\S\ref{sec:cutoffs} and Remark \ref{rem:cutoffs} below.

Let $\Hilbert_\ell$ be a closed subspace of the restriction of $\Hilbert$ to $\Omega_\ell$ (in particular, $\Hilbert_\ell$ may contain zero Dirichlet boundary conditions on part of $\partial\Omega_\ell$, in addition to any imposed via $\Hilbert$).
Let
\beq\label{eq:sesquiell_intro}
a_\ell (u,v):= \int_{\Omega_\ell} k^{-2} (\coeffA_\ell \nabla u)\cdot\overline{\nabla v} + k^{-1} (\coeffB_\ell\cdot\nabla u)\overline{v} - \coeffc_\ell^{-2} u\overline{v} - \ri k^{-1} \int_{\partial\Omega_\ell} \theta_\ell u\overline{v}.
\eeq
We assume that 
\beqs
\coeffA_\ell \equiv \coeffA, \quad \coeffB_\ell \equiv \coeffB, \quad \coeffc_\ell \equiv \coeffc,  \quad\tand \theta\equiv \theta_\ell \quad\ton \supp \chiellg,
\eeqs
so that $a(\cdot,\cdot) \equiv a_\ell(\cdot,\cdot)$ on $\supp \chiellg$, and, in particular, near $\partial\Omega\cap \partial\Omega_\ell$. This means that the local problems have the same boundary conditions on $\partial\Omega\cap \partial\Omega_\ell$ as the global problem, but away from 
$\supp \chiellg$ (i.e., near the boundaries of the subdomains that are not boundaries of $\Omega$) the local problems can have different boundary conditions and/or absorbing layers to the global problem.

We allow the coefficients $\coeffA, \coeffA_\ell, \coeffB, \coeffB_\ell, \coeffc,$ and $\coeffc_\ell$ to depend on $k$ (e.g., for a Cartesian PML near $\partial\Omega$, $\coeffB\sim k^{-1}$), but we assume that $\coeffA, \coeffA_\ell, \coeffB, \coeffB_\ell$ are all bounded above independent of $k$, and $\coeffA,\coeffA_\ell, \coeffc$ and $\coeffc_\ell$ are all bounded below independent of $k$.

Let 
\beq\label{eq:1knorm}
\N{u}^2_{H^1_k(\Omega)}:= k^{-2} \big\|\nabla u\big\|_{L^2(\Omega)}^2 + \N{u}^2_{L^2(\Omega)}.
\eeq
We note that many papers on the numerical analysis of the Helmholtz equation use the weighted $H^1$ norm $\vertiii{v}^2:=\N{\nabla v}^2_{L^2(\Omega)} + k^2\N{v}^2_{L^2(\Omega)}$; we work with \eqref{eq:1knorm} instead, because weighting the $j$th derivative with $k^{-j}$ makes the norm dimensionless and is easier to keep track of than weighting the $j$th derivative with $k^{-j+1}$.

The assumptions above imply that $a(\cdot,\cdot)$ is continuous and satisfies a G\aa rding inequality on $\Hilbert$ (equipped with the norm \eqref{eq:1knorm}), with the constants in these inequalities independent of $k$; similarly, $a_\ell(\cdot,\cdot)$ is continuous and satisfies a G\aa rding inequality on $\Hilbert_\ell$, again with constants independent of $k$.

\subsubsection{The discretised problem and preconditioners}

Let $\fine \subset \cV$ (the \emph{fine space}) be piecewise-polynomial Lagrange finite elements on a shape-regular mesh of size $h$, 
$\subdomain:= \fine|_{\Omega_\ell} \cap\Hilbert_\ell$, and let $\coarse$ (the \emph{coarse space}) be a subspace of $\fine$.

Let $\matrixA$ be the Galerkin matrix of $a(\cdot,\cdot)$ discretised in $\fine$; we assume that $h$ depends on the polynomial degree $p$ and $k$ in such a way that $\matrixA$ is invertible (we discuss these conditions on $h$ below Corollary \ref{cor:pwp}). 

Let $\matrixR_0$ 
be the restriction matrix from the degrees of freedom
on $\fine$ to the degrees of freedoms on $\coarse$, and 
 let $\matrixA_0:= \matrixR_0  \matrixA \matrixR_0 ^T$; i.e., 
$\matrixA_0$ is the Galerkin matrix of $a(\cdot,\cdot)$ discretised in $\cV_0$. 
Let $\matrixA_\ell$ be the Galerkin matrix of $a_\ell(\cdot,\cdot)$ discretised in $\cV_{\ell,h}$.

We define two 
restriction matrices from the degrees of freedom on $\fine$ to the degrees of freedoms on $\subdomain$, 
$\matrixR_\ell^{\chiell}$ and $\matrixR_\ell^{\chiellg}$, where 
$\matrixR_\ell^{\chiell}$ is weighted by $\chiell$, and $\matrixR_\ell^{\chiellg}$ is weighted by $\chiellg$ (see \eqref{eq:weightedR} below).
Let 
\beq\label{eq:matrixB}
\matrixB_L^{-1}= \matrixB_L^{-1}(\matrixA):= 
\matrixR_0^T\matrixAzero^{-1} \matrixR_0 
+
\bigg(\sum_{\ell=1}^N
\big(\matrixR_\ell^{\chiell}\big)^T\matrixA_\ell^{-1} \matrixR_\ell^{\chiellg}\bigg)
\big(\matrixI -\matrixA\matrixR_0^T\matrixAzero^{-1} \matrixR_0 \big).
\eeq
%
%
%
We call $\matrixB_L^{-1}$ a \emph{hybrid} Schwarz preconditioner because the coarse and subdomain solves are combined in a multiplicative way; this idea was first introduced in 
\cite{Ma:93, MaBr:96}.

Let the real symmetric positive-definite matrix $\matrixD_k$ 
be such that, for all $v_h\in\fine$ with coefficient vectors $\bV$,
\begin{equation}\label{eq:ip2}
\|v_h\|_{H^1_k(\Omega)}^2  = \langle \matrixD_k \bV,  \bV\rangle_2,
\eeq
where $\langle\cdot,\cdot\rangle_2$ denotes the Euclidean inner product.
Let $^\dagger$ denote the adjoint with respect to the Euclidean inner product $\langle\cdot,\cdot\rangle$; i.e., $\matrixA^\dagger= \overline{\matrixA}^T$. 
A few lines of calculation then show that 
\begin{gather}\label{eq:adjoint1}
\tif\quad \bW_j  =\matrixD_k \bV_j\quad \text{ then }\quad\big\langle \bV_1, \bV_2\big\rangle_{\matrixD_k}=\big\langle \bW_1, \bW_2\big\rangle_{\matrixD_k^{-1}}\quad\tand \\
\label{eq:adjoint2}
\big\langle \bV_1, \matrixB_L^{-1}(\matrixA)\matrixA \bV_2 \big\rangle_{\matrixD_k}
=\big\langle \matrixA^\dagger (\matrixB_L^{-1}(\matrixA))^\dagger \bW_1, \bW_2 \big\rangle_{\matrixD_k^{-1}}.
\end{gather}
We therefore set
$
\matrixB_R^{-1}(\matrixA):=(\matrixB_L^{-1}(\matrixA^\dagger))^\dagger,
$
i.e., 
\beqs
\matrixB_R^{-1}(\matrixA)
=\matrixR_0^T\matrixAzero^{-1} \matrixR_0 
+
\big(\matrixI -\matrixR_0^T\matrixAzero^{-1} \matrixR_0\matrixA\big)
\bigg(\sum_{\ell=1}^N
\big(\matrixR_\ell^{\chiellg}\big)^T\matrixA_\ell^{-1} \matrixR_\ell^{\chiell} \bigg),
\eeqs
so that, by \eqref{eq:adjoint1} and \eqref{eq:adjoint2} with $\matrixA$ replaced by $\matrixA^\dagger$,
\beq\label{eq:adjoint2a}
\big\langle \bV_1, \big(\matrixI-\matrixB_L^{-1}(\matrixA^\dagger)\matrixA^\dagger\big) \bV_2 \big\rangle_{\matrixD_k}
=\big\langle\big(\matrixI- \matrixA (\matrixB_R^{-1}(\matrixA))\big)\bW_1, \bW_2 \big\rangle_{\matrixD_k^{-1}}.
\eeq

\subsubsection{Informal statement of the main result}

\begin{theorem}\label{thm:informal1}
With the set-up above, suppose additionally that

(a) The subdomain diameters are all proportional to $k^{-1}$ and the subdomains all have generous overlap (i.e., the overlaps are also proportional to $k^{-1}$).

(b) The coarse space is such that the following holds:~if the Helmholtz problem is solved using the Galerkin method in the coarse space, 
    then 
    \bit
    \item the $H^1_k$ Galerkin error is bounded (independent of $k$) by the $H^1_k$ norm of the solution, and 
    \item the $L^2$ Galerkin error is bounded by a sufficiently-small (independent of $k$) multiple of the  $H^1_k$ norm of the solution. 
    \eit
    (We recall below that both these bounds hold if the Galerkin solution in the coarse space is proved to be quasi-optimal via the Schatz argument.) 
    
    Then
\beqs
\N{\matrixI-\matrixB_L^{-1}\matrixA  }_{\matrixD_k} \quad\text{ is small  (independent of $k$)}
\eeqs
and thus there exists $0<c<1$ (independent of $k$) such that 
\bit
\item 
the preconditioned fixed point iteration $\bx^{n+1} = \bx^{n} + \matrixB_L^{-1}({\bf b}-\matrixA \bx^n)$ for solving $\matrixA\bx={\bf b}$ satisfies
\beq\label{eq:plant1}
\N{\bx-\bx^n}_{\matrixD_k} \leq c^n \big\| \bx- \bx^0
\big\|_{\matrixD_k},
\eeq
\item when GMRES is applied in the $\matrixD_k$ inner product, the residual $\br^n:= {\bf b}- \matrixA \bx^n$ satisfies 
\beq\label{eq:plant2}
\N{\br^n}_{\matrixD_k} /\N{\br^0}_{\matrixD_k} \leq c^n,\qquad \tand
\eeq
\item when GMRES is applied in the Euclidean inner product and the fine mesh is quasi-uniform, the residual $\br^n$ satisfies 
\beq\label{eq:plant3}
\N{\br^n}_{2} /\N{\br^0}_{2} \leq c^n (hk)^{-1}.
\eeq
\eit
Furthermore, if Assumption (b) also holds for the adjoint sesquilinear form, then, by \eqref{eq:adjoint2a}, 
$\N{\matrixI-\matrixA\matrixB_R^{-1}}_{\matrixD_k^{-1}}$ is small, and so an analogous result holds for right preconditioning (in the $\matrixD_k^{-1}$ inner product). 

\label{thm:informal}
\end{theorem}

We make the following remarks. 

(i) Since the Helmholtz solution operator at high frequency involves propagation at length scales independent of $k$,  
and subdomains of size $k^{-1}$ cannot see this, the coarse space must resolve this propagation. 
Theorem \ref{thm:informal1} encodes this requirement in Assumption (b). We discuss below -- in the context of piecewise polynomials -- how one might seek to weaken Assumption (b), but highlight 
that the coarse spaces of  \cite{HuLi:24, LuXuZhZo:24, MaAlSc:24, FuGoLiWa:24} all satisfy Assumption (b) (see \S\ref{sec:coarse_other} below).

(ii) Following on from (i), once the wave nature of the solution is resolved by the coarse space, one obtains the Galerkin solution at the accuracy of the fine space.
One can therefore view the role of the subdomains as inexpensive high-accuracy interpolators, tailored to the specific wave problem.
%
%

(iii) The main ways in which Theorem \ref{thm:informal1} is a simplification of the main result (Theorem \ref{thm:main1} below) are that 
\bit
\item Theorem \ref{thm:main1} allows arbitrary-sized subdomains. We have considered subdomains of size $k^{-1}$ in Theorem \ref{thm:informal1} since, in practice, one wants to take the subdomains small for parallel scalability. Indeed, with $h$ chosen as a function of $k$ and $p$ to maintain accuracy as $k\to\infty$, the number of degrees of freedom in each $k^{-1}$-sized subdomain grows mildly with $k$ (like $k^{d/p}$ or $k^{d/(2p)}$ for nontrapping problems, depending on the measure of accuracy used).

\item Theorem \ref{thm:informal1} assumes that the fine space consists of piecewise polynomials; the only assumption on the fine space in Theorem \ref{thm:main1} is a super-approximation result (Assumption \ref{ass:super} below) which is satisfied by piecewise polynomials, but also, in principle, other spaces.
\eit


(iv) The result \eqref{eq:plant3} in Theorem \ref{thm:informal1} (about standard GMRES) follows from the result \eqref{eq:plant2} about weighted GMRES  via an inverse estimate, as recognised in \cite[Corollary 5.8]{GoGrSp:20}.
The numerical experiments in 
\cite[Experiment 1]{GSV1}, \cite[\S6]{BDGST} showed little difference in the number of weighted/unweighted GMRES iterations.

(v) The set up of Theorem \ref{thm:informal1} allows a wide variety of boundary conditions/absorbing layers on the global and subdomain problems (including, e.g., a PML). The existing two-level hybrid Schwarz analyses in \cite{HuLi:24, LuXuZhZo:24, MaAlSc:24, FuGoLiWa:24} all consider impedance boundary conditions on $\partial\Omega$ (which give $k$-independent errors when approximating the Sommerfeld radiation condition \cite{GLS1}) and either impedance or Dirichlet boundary conditions on $\partial\Omega_\ell$.


(vi) By \eqref{eq:matrixB}, 
\begin{align*}
\matrixI-\matrixB_L^{-1}\matrixA&=
(\matrixI-\matrixR_0^T\matrixAzero^{-1} \matrixR_0 \matrixA)
-
\bigg(\sum_{\ell=1}^N
\big(\matrixR_\ell^{\chiell}\big)^T\matrixA_\ell^{-1} \matrixR_\ell^{\chiellg}\matrixA\bigg)
\big(\matrixI -\matrixR_0^T\matrixAzero^{-1} \matrixR_0 \matrixA\big)\\
&=
\bigg(\matrixI- \sum_{\ell=1}^N
\big(\matrixR_\ell^{\chiell}\big)^T\matrixA_\ell^{-1} \matrixR_\ell^{\chiellg}\matrixA\bigg)
\big(\matrixI -\matrixR_0^T\matrixAzero^{-1} \matrixR_0 \matrixA\big).
\end{align*}
One way to understand heuristically why $\matrixI-\matrixB_L^{-1}\matrixA$ is small is the following 
\bit
\item $\matrixI -\matrixR_0^T\matrixAzero^{-1} \matrixR_0 \matrixA$ 
is effective at removing frequencies $\lesssim k$, i.e., it acts like a high-pass filter, and 
\item $\matrixI- \sum_{\ell=1}^N
\big(\matrixR_\ell^{\chiell}\big)^T\matrixA_\ell^{-1} \matrixR_\ell^{\chiellg}\matrixA$ is effective at removing frequencies $\gtrsim k$, i.e., it acts like a low-pass filter;
\eit
thus the product of these two operators is small. 

Regarding the first point:~this is because the Galerkin solution in the coarse space is quasi-optimal and the best approximation of a frequency $\lesssim k$ function is small \cite{G1}.
 
Regarding the second point:~at the continuous level, frequencies $\gg k$ do not propagate under the action of the Helmholtz solution operator (as a consequence of semiclassical ellipticity \cite[Theorem E.33]{DyZw:19}). Thus, modulo frequencies $\lesssim k$, 
the solution can be well-approximated by a sum of local solves. Provided the local spaces $\Hilbert_{\ell,h}$ accurately approximate $\Hilbert_h$ restricted to $\Omega_\ell$, this property is inherited by the discrete solution operator. (We encode this accurate-approximation assumption in the super-approximation-type estimate, Assumption~\ref{ass:super} below.)

\subsubsection{The main result applied to piecewise-polynomial coarse spaces}

Let $\Csol$ denote the $L^2\to H^1_k$ norm of the Helmholtz solution operator (recall that this $\sim kL$ when the problem is nontrapping \cite[\S4.6]{DyZw:19}, where $L$ is the characteristic length scale of $\Omega$).

\begin{corollary}[Piecewise-polynomial coarse space of fixed degree]\label{cor:pwp}
Suppose that $\fine$ consists of degree-$p$ piecewise polynomials on a mesh of size $h$ and $\coarse\subset\fine$ with the coarse mesh size $\Hcoarse$ and 
each element of the coarse mesh a union of elements of the fine mesh. 

Suppose that the subdomain diameters are all proportional to $k^{-1}$ and the subdomains all have generous overlap (i.e., the overlaps are also proportional to $k^{-1}$).

If the domain and coefficients have suitable regularity, 
\beq\label{eq:qo_threshold}
(k\Hcoarse)^{p}\Csol \text{ is sufficiently small},
\eeq
and, 
given a coarsening factor $\Ccoarse>1$,
\beqs
h=\Hcoarse/\Ccoarse, 
\eeqs
then the properties \eqref{eq:plant1}, \eqref{eq:plant2}, and \eqref{eq:plant3} hold.
\end{corollary}

The precise statement of Corollary \ref{cor:pwp} is Corollary \ref{cor:pwp2} below. 
We make the following remarks about Corollary \ref{cor:pwp}/Corollary \ref{cor:pwp2}.

(i) The condition \eqref{eq:qo_threshold} on $\Hcoarse$ means that the sequence of Galerkin solutions in the sequence of coarse spaces is quasi-optimal. 
One could hope to prove a result with the condition \eqref{eq:qo_threshold} replaced by ``$(k\Hcoarse)^{2p}\Csol$ sufficiently small", with this threshold ensuring that, for $k$-oscillatory data, the Galerkin solution in the coarse space has bounded relative error  \cite{GS3}. However, as discussed after Theorem \ref{thm:informal1}, the coarse problem needs to sufficiently resolve the wave-nature of the solution. 
We note that the experiments in 
\cite[Table 9]{BoDoJoTo:21} show that, when a two-level additive Schwarz preconditioner with a fixed number of points per wavelength in both fine and coarse space is used, the number of GMRES iterations is large if the coarse mesh does not contain a sufficient number of points per wavelength.

(ii) In the set up of Corollary \ref{cor:pwp}, 
\beqs
\frac{\text{coarse-space dimension }
}{\text{fine-space dimension }
} \sim  
\bigg(\frac{h }{\Hcoarse}
\bigg)^d
\sim \bigg(\frac{1}{\Ccoarse}
\bigg)^d
\tas k\to \infty.
\eeqs
The coarse problem is still large -- as expected from the requirement discussed in (i) that the coarse problem resolve the oscillations -- and, in practice, is often solved using a one-level DD method  (see, e.g., \cite{ToJoDoHoOpRi:22}, \cite[\S6]{BDGST}). 
We do not attempt to analyse this set up in the present paper (in common with the two-level analyses in 
 \cite{HuLi:24, LuXuZhZo:24, MaAlSc:24, FuGoLiWa:24}).

 (iii) As discussed after Theorem \ref{thm:informal1}, once the wave nature of the solution is resolved by the coarse space, one obtains the Galerkin solution at the accuracy of the fine space.
The concrete realisation of this in the setting of Corollary \ref{cor:pwp} is that, for fixed large $k$, once $\Hcoarse$ 
satisfies \eqref{eq:qo_threshold}, arbitrary accuracy can be achieved by increasing $\Ccoarse$ (so that $h/\Hcoarse\to 0$).

(iv) The follow-up paper \cite{IvanEuan1} proves results about piecewise-polynomial coarse spaces with polynomial degree increasing like $\log k$. This analysis is restricted to a specific Helmholtz problem (where the radiation condition is approximated by a complex absorbing potential \cite{RiMe:93}), no partition of unity functions appear in the preconditioner, and the subdomains have zero Dirichlet boundary conditions and width a sufficiently small multiple of $k^{-1}$ (so that the subdomain problems are coercive).
The difficulty in applying the theory in the present paper with increasing polynomial degree is the super-approximation result in Assumption \ref{ass:super}/Lemma \ref{lem:super}; the current proof of this result gives a constant that blows up as the polynomial degree increases.

\subsection{Discussion about the cutoffs $\chiell$ and $\chiellg$}\label{sec:cutoffs}

The one-level part of the preconditioner \eqref{eq:matrixB} is
$\sum_{\ell=1}^N
\big(\matrixR_\ell^{\chiell}\big)^T\matrixA_\ell^{-1} \matrixR_\ell^{\chiellg}$

In the literature, one-level methods of the form 
\begin{subequations} \label{eq:group}
\noindent
\begin{minipage}{0.45\linewidth}
\begin{equation}
\sum_{\ell=1}^N
\big(\matrixR_\ell^{\chiell}\big)^T\matrixA_\ell^{-1} \matrixR_\ell^{\chiell}, \quad\tor
 \tag{\ref{eq:group}a} \label{eq:group:a}
\end{equation}
\end{minipage}
\hfill
\begin{minipage}{0.45\linewidth}
\begin{equation}
\sum_{\ell=1}^N
\big(\matrixR_\ell^{\chiell}\big)^T\matrixA_\ell^{-1} \matrixR_\ell 
\tag{\ref{eq:group}b} \label{eq:group:b}
\end{equation}
\end{minipage}
\end{subequations}
often appear, where $\matrixR_\ell$ is just restriction to $\Omega_\ell$.
For example, \eqref{eq:group:a} appears in the one-level analysis of \cite{GSZ1} (see \cite[Equation 1.13]{GSZ1}) and as the one-level component of the two-level methods in \cite{MaAlSc:24} and \cite{FuGoLiWa:24} (see \cite[Equation 3.9]{MaAlSc:24} and \cite[Page 7]{FuGoLiWa:24}).
Furthermore \eqref{eq:group:b} appears in the ``ORAS'' one-level preconditioner (see, e.g., \cite[Equation 10]{BoDoJoTo:21} and \cite[Equation 2.5]{GoGrSp:23}) and as the one-level component of the two-level methods in \cite{HuLi:24} and \cite{LuXuZhZo:24} (see \cite[Equation 3.4]{HuLi:24} and \cite[Equations 4.30 and 4.31]{LuXuZhZo:24}).

Comparing \eqref{eq:matrixB} and \eqref{eq:group}, we see that the one-level part preconditioner \eqref{eq:matrixB} in the present paper can be seen as somewhere ``in between" \eqref{eq:group:a} and \eqref{eq:group:b} -- the restriction operator at the end of $\sum_{\ell=1}^N
\big(\matrixR_\ell^{\chiell}\big)^T\matrixA_\ell^{-1} \matrixR_\ell^{\chiellg}$
 acts on a large set that ${\rm supp} \chiell$, but not all of $\Omega_\ell$.

The technical reason why the proof of the main result does not hold if the one-level part is either \eqref{eq:group:a} or \eqref{eq:group:b} is described in Remark \ref{rem:cutoffs}.
We therefore see the presence of $\chiellg$ as (informally) ``the price we pay'' to have this simple theory of two-level methods -- requiring only that the coarse solutions are quasi-optimal, but not relying on any specific details of the coarse space, and holding for arbitrary absorbing
layers/boundary conditions on both the global and local Helmholtz problems.

\paragraph{Outline of the paper}
\S\ref{sec:ass} states our abstract assumptions 
and \S\ref{sec:3} states and proves the main result (i.e., the rigorous version of Theorem \ref{thm:informal1}) under these assumptions. 
\S\ref{sec:verify_ass} shows how the abstract assumptions -- excluding that on the coarse space -- are valid for Helmholtz problems with piecewise-polynomial fine spaces, and 
\S\ref{sec:coarse} then gives conditions under which the assumption on the coarse space  holds.
\S\ref{sec:6} gives the rigorous statement of  Corollary \ref{cor:pwp}.
\S\ref{sec:appendix} gives the matrix form of the preconditioner, and \S\ref{sec:Schatz_app} recaps the Schatz argument and Aubin--Nitsche lemma.

\section{Statement of the abstract assumptions}\label{sec:ass}

\subsection{The space $\Hilbert$ and sesquilinear form $a(\cdot,\cdot)$}

Let $\Omega\subset\Rea^d$ be a bounded Lipschitz domain, and let $\Hilbert$ be a closed subspace of $H^1(\Omega)$.\label{sec:ass1}

Let $a:\Hilbert\times\Hilbert\to \Com$ be a continuous sesquilinear form; i.e., given $k_0>0$ there exists $\Ccont>0$ such that for all $k\geq k_0$
\beqs
\big|a(u,v)\big| \leq \Ccont \N{u}_{H^1_k(\Omega)}\N{v}_{H^1_k(\Omega)} \quad\tfa u,v\in \Hilbert.
\eeqs
Assume that, given $F\in \Hilbert^*$, there exists a unique solution to the variational problem 
\beq\label{eq:vp}
\tfind u\in \Hilbert \,\,\tst \,\, a(u,v) =F(v)\,\,\tfa v\in \Hilbert.
\eeq

\subsection{The space $\Hilbert_h$, subdomains $\Omega_\ell$, and spaces $\Hilbert_\ell$ and $\Hilbert_{\ell,h}$}\label{sec:ass2}

For each $k>0$, let $\{ \Omega_\ell(k) \}^N_{\ell =1}$ (the \emph{subdomains}) form an overlapping cover of $\Omega$, with each $\Omega_\ell(k)$ a non-empty open Lipschitz domain with characteristic length scale $H_\ell(k)$ (i.e., its diameter $\sim  H_\ell(k)$, its surface area $\sim  H_\ell(k)^{d-1}$, and its volume $\sim H_\ell(k)^d$). For brevity, we omit the dependence on $k$ in this notation and write, e.g., $\Omega_\ell$ instead of $\Omega_\ell(k)$, but we emphasise that the subdomains can depend on $k$.
Let 
\beq\label{eq:Lambda}
\Lambda  := \max\bigl\{ \#\Lambda (\ell ):\ell =1,...,N\bigr\},\quad\text{where}\quad \Lambda (\ell ):=\bigl\{ \ell^\prime  :\Omega_\ell \cap \Omega_{\ell^\prime}  \not =\emptyset \bigr\};
\eeq
i.e., $\Lambda$ is the maximum number of subdomains that can overlap any given subdomain.

Define the norm $\|\cdot\|_{H^1_k(\Omega_\ell)}$ and associated inner product in an analogous way to \eqref{eq:1knorm} with $\Omega$ replaced by $\Omega_\ell$. 
Let $\cR_\ell: L^2(\Omega)\to L^2(\Omega_\ell)$ denote the restriction operator, and let $\cR_\ell^*$ denote its $L^2$ adjoint; i.e., extension by zero. Note that 
also $\cR_\ell: H^1(\Omega)\to H^1(\Omega_\ell)$.

Let $\Hilbert_\ell$ be a closed subspace of $\cR_\ell(\Hilbert)$ and 
let $\{\Hilbert_{\ell,h}\}_{h>0}$ be a family of finite-dimensional subspaces of $\Hilbert_\ell$, indexed by $h>0$.
Let $\{\chiell\}_{\ell=1}^N$ be a partition of unity subordinate to $\{ \Omega_\ell \}^N_{\ell =1}$ (so that, in particular, $\cR_\ell^* \chi_\ell \cR_\ell = \chi_\ell$). Let $\{\chiellg\}_{\ell=1}^N$ be such that
for $\ell=1,\ldots,N$, 

(i) $\chiellg \in L^\infty(\Omega_\ell;[0,1])$, 

(ii) $\chiellg \equiv 1$ on $\supp \chiell$, 

(iii) $\chiell, \chiellg : \cR_\ell (\Hilbert) \to \Hilbert_\ell \cap \Hot(\Omega_\ell)$, where 
\beqs
\Hot(\Omega_\ell) := \big\{ v \in H^1(\Omega_\ell) \text{ with } v= 0 \text{ on } \partial\Omega_\ell\setminus\partial\Omega\big\}.
\eeqs

Let $\{\Hilbert_h\}_{h>0}$ be a family of finite-dimensional subspaces of $\Hilbert$, indexed by $h>0$.

\begin{assumption}[Super-approximation-type assumption]\label{ass:super}
There exist linear operators $\cI_{h}^{\ell}$ such that  $\cI_h^\ell \chiell,\cI_h^\ell \chiellg: (\cR_\ell \fine + \Hilbert_{\ell,h}) \to \Hilbert_{\ell,h} \cap \Hot(\Omega_\ell)$ 
and there exist $\mu_\ell=\mu_\ell(k):(0,\infty)\to (0,\infty)$, $\ell=1,\ldots,N$, such that for all $v\in (\cR_\ell \fine + \Hilbert_{\ell,h})$ and $k\geq k_0$
\begin{align}
\label{eq:super_approx}
\max\Big\{\big\| (I-\cI_{h}^{\ell})(\chiell v)\big\|_{H^1_k(\Omega_\ell)},\,
\big\| (I-\cI_{h}^{\ell})(\chiellg v)\big\|_{H^1_k(\Omega_\ell)}\Big\}
\leq \mu_\ell\N{v}_{H^1_k(\Omega_\ell)}.
\end{align}
\end{assumption}

We make Assumption \ref{ass:super} because, in the proof of Theorem \ref{thm:main1}, $\cI_h^\ell \chiell$ and $\cI_h^\ell \chiellg$ are applied \emph{both} to functions in $\mathcal{R}_\ell \mathcal{V}_h$
\emph{and} to functions in $\mathcal{V}_{\ell,h}$. In many situations, these spaces are not the same because of boundary conditions:~e.g., with a PML/absorbing layer on the subdomains, $\mathcal{V}_{\ell,h}$ has a zero boundary condition on $\partial\Omega_\ell$, but $\mathcal{R}_\ell \mathcal{V}_h$ has no boundary condition on $\partial\Omega_\ell\setminus \partial\Omega$.

Lemma \ref{lem:super}
below shows that \eqref{eq:super_approx} holds when $\Hilbert_h$ consists of piecewise polynomials on a shape-regular mesh; in this case $\mu_\ell:= h_\ell/\delta_\ell$, where $h_\ell  := \max_{\element   \subset \overline{\Omega_\ell}}  h_\element $ and $\delta_\ell$ is such that $\|\partial^\alpha \chiell\|_{L^\infty(K)}, \|\partial^\alpha \chiellg\|_{L^\infty(K)} \leq C \delta_\ell^{-|\alpha|}$ (and thus $\delta_\ell$ is related to the overlap of the subdomains).
Let 
\beq\label{eq:mu}
\mu:= \max_{\ell} \mu_\ell.
\eeq

\subsection{The local sesquilinear forms $a_\ell(\cdot,\cdot)$}\label{sec:ass3}

Let $a_\ell:\cR_\ell(\Hilbert)\times\cR_\ell(\Hilbert)\to \Com$ be a continuous sesquilinear form depending on $k$. 
Let $\delta_\ell= \delta_\ell(k): (0,\infty)\to (0,1), \ell=1\ldots N$ 
and $\gamma_\ell= \gamma_\ell(k): (0,\infty)\to (0,1), \ell=1\ldots N$. 
Let 
\beq\label{eq:delta}
\delta:=\min_{\ell} \delta_\ell,
\quad\tand\quad \gamma:= \max_{\ell} \gamma_\ell.
\eeq

We assume that the continuity constant of $a_\ell(\cdot,\cdot)$ is the same as for $a(\cdot,\cdot)$; i.e., for all $k\geq k_0$
\beqs
\big|a_\ell(u,v)\big| \leq \Ccont \N{u}_{H^1_k(\Omega_\ell)}\N{v}_{H^1_k(\Omega_\ell)} \quad\tfa u,v\in \cR_\ell(\Hilbert).
\eeqs

\begin{assumption}\label{ass:dis}
$a_\ell(\cdot,\cdot)$ satisfies the discrete inf-sup condition
\beq\label{eq:is}
\inf_{u_{\ell,h}\in \Hilbert_{\ell,h}} \sup_{v_{\ell,h}\in\Hilbert_{\ell,h}}\frac{
\big| a_\ell(u_{\ell,h}, v_{\ell,h})\big|
}{
\N{u_{\ell,h}}_{H^1_k(\Omega_\ell)} 
\N{v_{\ell,h}}_{H^1_k(\Omega_\ell)} 
}
\geq \gamma_\ell^{-1} \quad\tfa k\geq k_0.
\eeq
\end{assumption}

\begin{assumption}[$a$ and $a_\ell$ agree ``in the interior of $\Omega_\ell$ and near $\partial\Omega$"]\label{ass:supports}
\beq\label{eq:supports}
a_\ell\big( \cR_\ell u_h, \cI_{h}^{\ell} (\chiellg v_{\ell,h})\big) = a\big(u_h, \cR_\ell^*\cI_{h}^{\ell}(\chiellg v_{\ell,h})\big)\quad \tfa u_h \in \Hilbert_h, v_{\ell,h}\in \Hilbert_{\ell,h}.
\eeq
\end{assumption}

\begin{assumption}[Bound on the derivative of $\chiell$]
There exists $C_{\rm PoU}>0$ such that for all $k\geq k_0$
\beq\label{eq:PoU_new1}
\N{\nabla\chiell}_{L^\infty(\element )}
\leq C_{\rm PoU} \delta_\ell^{-1}.
\eeq
\end{assumption}

\begin{assumption}[Commutator bound]\label{ass:commutator}
Given $k_0>0$ there exists $\Ccom>0$ 
such that for all $k\geq k_0$, all $u,v\in \cR_\ell(\Hilbert)$, and $\ell=1,\ldots,N$,
\begin{align}\label{eq:commutator}
\big| a_\ell(\chiellg u, v)- a_\ell(u,\chiellg v)\big|\leq \Ccom (k\delta_\ell)^{-1}\big(1+(k\delta_\ell)^{-1}\big)\N{u}_{L^2(\Omega_\ell)} \N{v}_{H^1_k(\Omega_\ell)}.
\end{align}
\end{assumption}

\subsection{The coarse space $\coarse$}\label{sec:ass_coarse}

\begin{assumption}\label{ass:coarse}
$\coarse$ is a finite-dimensional subspace of $\Hilbert_h$ and 
there exists an operator $Q_0: \fine\to \coarse$ and $\sigma_{L^2}, \sigma_{H^1}:(0,\infty)\to (0,\infty)$ such that for all $k\geq k_0$
\beq\label{eq:sigma}
\N{(I-Q_0)v_h}_{L^2(\Omega)}\leq \sigmat(k) \N{v_h}_{H^1_k(\Omega)},\quad
\N{(I-Q_0)v_h}_{H^1_k(\Omega)}\leq \sigmao(k) \N{v_h}_{H^1_k(\Omega)}.
\eeq
\end{assumption}

We make two remarks about Assumption \ref{ass:coarse}. First, one way to define $Q_0$ is as the Galerkin solution in $\coarse$:
\beq\label{eq:Q0}
a\big( Q_0 u_h,v_0 \big)= a\big( u_h,v_0 \big)\quad\tfa v_0 \in \coarse,
\eeq
but the main abstract result (Theorem \ref{thm:main1}) only requires Assumption \ref{ass:coarse} (and thus covers, e.g., the LOD coarse space of \cite{LuXuZhZo:24}, where the projection $Q_0$ is defined via a Petrov--Galerkin method). Second, for the result of the main abstract result (Theorem \ref{thm:main1} below) to be useful, $\sigmat$ in \eqref{eq:sigma} must be small.

\section{The main abstract result and its proof}\label{sec:3}

\subsection{Statement of the main abstract result}

For $\ell=1,\ldots, N$, let $Q_\ell: \Hilbert_h\to \Hilbert_{\ell,h}$ be defined by 
\beq\label{eq:Qell}
a_{\ell}\big( Q_\ell u_h, v_{\ell,h} \big)= a\big( u_h, \cR_\ell^*\cI_{h}^{\ell}(\chiellg v_{\ell,h}) \big)\quad\tfa v_{\ell,h} \in \subdomain.
\eeq
Let 
\beq\label{eq:Q}
Q=Q_0 + \sum_{\ell=1}^N \cR_\ell^*\cI_{h}^{\ell}\big(\chi_\ell Q_\ell(I-Q_0)\big),
\eeq
where $Q_0$ is as in Assumption \ref{ass:coarse}.

\S\ref{sec:appendix} below shows that when $Q_0$ is defined by \eqref{eq:Q0} (i.e., the Galerkin solution in the coarse space) then the matrix form of $Q$ is $\matrixB_L^{-1}\matrixA$ with $\matrixB_L^{-1}$ defined by \eqref{eq:matrixB}.

\begin{theorem}[Bound on $I-Q$]\label{thm:main1}
Under the assumptions in \S\ref{sec:ass} and with $Q$ defined by \eqref{eq:Q}, for all $v_h\in\Hilbert_h$ and $k\geq k_0$
\beq\label{eq:main}
\begin{aligned}
\N{(I-Q)v_h}_{H^1_k(\Omega)} &\leq 2\Lambda \Big[(1+\mu) 
\big(1 + C_{\rm PoU}(k\delta)^{-1})\big)
\gamma 
\Ccom (k\delta)^{-1} \big(1+ (k\delta)^{-1}\big) \sigmat \\
&\hspace{0.25cm}+\big(1+(1+\mu)\big(1 + C_{\rm PoU}(k\delta)^{-1})\big) (1+\gamma  \Ccont) \big)\mu \sigmao
\Big]\N{v_h}_{H^1_k(\Omega)}.
\end{aligned}
\eeq
\end{theorem}

The key point is that the right-hand side of \eqref{eq:main} can be made small by making both $\sigmat$ and $\mu$ small.

\bre[The relationships between $\coarse, \fine$, and $\Hilbert_{\ell,h}$]
In conventional two-level DD theory
\beqs
\fine= \sum_{\ell=1}^N \cR^*_\ell \Hilbert_{\ell,h}, 
\eeqs
since, modulo boundary conditions on $\partial\Omega_\ell\setminus \partial\Omega$, $\cR_\ell \coarse\subset \Hilbert_{\ell,h}$, and so one can take a single interpolation operator $\cI_h$ in Assumption \ref{ass:super}, rather than one for each subdomain.
However, the assumptions above allow for 
\beqs
\fine= \coarse + \sum_{\ell=1}^N \cR^*_\ell \Hilbert_{\ell,h};
\eeqs
that is, the local spaces need not \emph{contain} the restrictions to $\Omega_\ell$ of the coarse space. Instead the local spaces must \emph{accurately approximate} the restrictions of the coarse space, via Assumption \ref{ass:super}. For instance, if $\Hilbert_0$ consists only of  functions oscillating at frequency $\lesssim k$ (i.e., such that $\|u\|_{H_k^{\ell}(\Omega)}\leq C_\ell \|u\|_{H_k^1(\Omega)}$, $\ell=2,\dots, p+1$) then one can take as the local space piecewise polynomials on a sufficiently small mesh. 
\ere

\subsection{Proof of Theorem \ref{thm:main1}}

\begin{lemma}[Consequences of the definition of $\Lambda$]
For all $v\in L^2(\Omega)$ and $w\in H^1(\Omega)$
\beq\label{eq:overlap2a}
\sum_{\ell=1}^N \N{v}^2_{L^2(\Omega_\ell)} \leq \Lambda \N{v}^2_{L^2(\Omega)}
\quad\tand\quad \sum_{\ell=1}^N \N{w}^2_{H^1_k(\Omega_\ell)} \leq \Lambda \N{w}^2_{H^1_k(\Omega)}.
\eeq
Furthermore, given $v_\ell\in \Hilbert_\ell$,
\beq\label{eq:overlap_new}
\bigg\|
\sum_{\ell=1}^N v_\ell 
\bigg\|^2_{H^1_k(\Omega)}
\leq 2\Lambda \sum_{\ell=1}^N \N{v_\ell}^2_{H^1_k(\Omega_\ell)}.
\eeq
\end{lemma}

\bpf
The bounds in \eqref{eq:overlap2a} and 
follow immediately from the definition \eqref{eq:Lambda} of $\Lambda$.
The bound \eqref{eq:overlap_new} without an explicit expression for the constant is proved in \cite[Lemma 4.2]{GSV1}. The definition of $\Lambda$ implies that \cite[Equation 4.8]{GSV1} holds with $\lesssim \sum_{\ell=1}^N \|v_\ell\|^2_{H^1_k(\Omega)}$ at the end replaced by $\leq \Lambda \sum_{\ell=1}^N \|v_\ell\|^2_{H^1_k(\Omega)}$. The result then follows, with the factor of $2$ arising from use of the inequality 
\beq\label{eq:quadratic}
(a+b)^2 \leq 2 (a^2 + b^2) \quad\tfa a,b>0,
\eeq
at the end of  \cite[Proof of Lemma 4.2]{GSV1} (with this constant hidden in the notation $\lesssim$ in \cite[Proof of Lemma 4.2]{GSV1}). 
\epf

\ble\label{lem:GSZ}
With $Q_\ell$ defined by \eqref{eq:Qell},
for all $u_h\in \Hilbert_h$,
\begin{align}\nonumber
&\N{Q_\ell u_h - 
\chiellg \cR_\ell u_h}_{H^1_k(\Omega_\ell)} \\
&\quad\leq\Big(\gamma_\ell \Ccom
(k\delta_\ell)^{-1} \big(1+ (k\delta_\ell)^{-1}\big) \N{u_h}_{L^2(\Omega_\ell)} + \big(1+\gamma_\ell\Ccont \big)\mu_\ell \N{u_h}_{H^1_k(\Omega_\ell)}
\Big).
\label{eq:GSZ}
\end{align}
\ele

\bre
The bound \eqref{eq:GSZ} is similar to that in \cite[Lemma 3.8]{GSZ1} except the latter contains only $\|u_h\|_{H^1_k(\Omega)}$ on the right-hand side.
\ere

\bpf[Proof of Lemma \ref{lem:GSZ}]
First observe that it is sufficient to prove that
\begin{align}\nonumber
&\N{Q_\ell u_h - \cI_{h}^{\ell}(\chiellg \cR_\ell u_h)}_{H^1_k(\Omega_\ell)} \\
&\qquad\leq \gamma_\ell \Big( \Ccom
(k\delta_\ell)^{-1} \big(1+ (k\delta_\ell)^{-1}\big) \N{u_h}_{L^2(\Omega_\ell)} + \Ccont \mu_\ell \N{u_h}_{H^1_k(\Omega_\ell)}
\Big),
\label{eq:GSZ2}
\end{align}
since then \eqref{eq:GSZ} follows from \eqref{eq:super_approx} and the triangle inequality.

By  the definition of $Q_\ell$ \eqref{eq:Qell} and the property \eqref{eq:supports},
\begin{align}\nonumber
&a_\ell \big(Q_\ell u_h - \cI_{h}^{\ell}(\chiellg\cR_\ell u_h), v_{\ell,h}\big) \\ 
&=a\big(u_h, \cR_\ell^*\cI_{h}^{\ell}(\chiellg v_{\ell,h})\big) - a_\ell\big( \cI_{h}^{\ell}(\chiellg \cR_\ell u_h), v_{\ell,h})
\label{eq:nomagic4}
\\ \nonumber
& =a_\ell\big(\cR_\ell u_h, \cI_{h}^{\ell}(\chiellg v_{\ell,h})\big) - a_\ell\big( \cI_{h}^{\ell}(\chiellg \cR_\ell u_h), v_{\ell,h})\\ \nonumber
&=a_\ell\big(\cR_\ell u_h, \chiellg v_{\ell,h}\big) - a_\ell\big( \chiellg \cR_\ell u_h, v_{\ell,h})
- a_\ell\big( \cR_\ell u_h, (I-\cI_{h}^{\ell})(\chiellg v_{\ell,h})\big) \\ \nonumber
&\hspace{6cm}
+ a_\ell\big( (I-\cI_{h}^{\ell})(\chiellg \cR_\ell u_h),  v_{\ell,h}\big) .
\end{align}
Then, by \eqref{eq:commutator} and \eqref{eq:super_approx},
\begin{align*}
&\big|a_\ell \big(Q_\ell u_h - \cI_{h}^{\ell}(\chiellg \cR_\ell u_h), v_{\ell,h}\big) \big|\leq 
\Ccom (k\delta_\ell)^{-1}\big(1+(k\delta_\ell)^{-1}\big)\N{u_h}_{L^2(\Omega_\ell)} \N{v_{\ell,h}}_{H^1_k(\Omega_\ell)} 
\\
&\hspace{6cm}+ \Ccont \mu_\ell \N{u_h}_{H^1_k(\Omega_\ell)} \N{v_{\ell,h}}_{H^1_k(\Omega_\ell)};
\end{align*}
the result then follows from \eqref{eq:is}.
\epf

\ble
\beq\label{eq:inter_bounded}
\N{\cI_{h}^{\ell}(\chiell \cR_\ell v_h)}_{H^1_k(\Omega_\ell)}\leq (1+\mu_\ell)\big(1 + C_{\rm PoU}(k\delta_\ell)^{-1}\big)\N{v_h}_{H^1_k(\Omega_\ell)} \quad\tfa v_h\in \Hilbert_h.
\eeq
\ele

\bpf
The result follows from \eqref{eq:super_approx},  the triangle inequality, the product rule, and \eqref{eq:PoU_new1}.
\epf

\bpf[Proof of Theorem \ref{thm:main1}]
By the facts that $\{\chiell\}_{\ell=1}^N$ is a partition of unity, $\cR_\ell^* \chiell \cR_\ell = \chiell$, and
 $\chiellg\equiv 1$ on $\supp \chiell$,  for all $u_h\in \Hilbert_h$,
\begin{align}\nonumber
u_h= \sum_{\ell=1}^N\cR_\ell^*\chi_{\ell}\cR_\ell u_h &= \sum_{\ell=1}^N\cR_\ell^*\Big(\cI_{h}^{\ell}(\chi_\ell \cR_\ell u_h)+
(I-\cI_{h}^{\ell})\chi_{\ell}\cR_\ell u_h\Big)
\\
&=\sum_{\ell=1}^N\cR_\ell^*\Big(\cI_{h}^{\ell}(\chi_\ell\chiellg \cR_\ell u_h)+
(I-\cI_{h}^{\ell})\chiell\cR_\ell u_h\Big). \label{eq:magic1}
\end{align}
By  \eqref{eq:Q} and the last displayed equation (with $u_h=(I-Q_0)v_h$), for $v_h\in \Hilbert_h$,
\begin{align*}
(Q- &I)v_h = -(I-Q_0)v_h+ \sum_{\ell=1}^N\cR_\ell^* \cI_{h}^{\ell}\big(\chi_\ell Q_\ell (I-Q_0)v_h\big)\\
&= \sum_{\ell=1}^N\cR_\ell^*\Big(- \cI_{h}^{\ell} (\chiell\chiellg \cR_\ell(I-Q_0)v_h)+ \cI_{h}^{\ell}\big(\chiell Q_\ell (I-Q_0)v_h\big)\\
&\hspace{6.3cm} -(I-\cI_{h}^{\ell})\chiell\cR_\ell(I-Q_0)v_h\Big)\\
&= \sum_{\ell=1}^N\cR_\ell^*\Big( \cI_{h}^{\ell}\big(\chiell \big(Q_\ell-\chiellg \cR_\ell\big)(I-Q_0)v_h\big)-(I-\cI_{h}^{\ell})\chiell\cR_\ell(I-Q_0)v_h\Big).
\end{align*}
Therefore, by (in this order) \eqref{eq:overlap_new}, \eqref{eq:inter_bounded}, \eqref{eq:super_approx}, 
 \eqref{eq:GSZ}, \eqref{eq:quadratic}, the definitions of $\mu$ \eqref{eq:mu} and $\delta$ \eqref{eq:delta}, and finally the two bounds in \eqref{eq:overlap2a}, 
\begin{align}\nonumber
\big\|&(Q- I)v_h \big\|^2_{H^1_k(\Omega)}\\ \nonumber
&\leq 2\Lambda \sum_{\ell=1}^N\Big[\big\| \cI_{h}^{\ell}\big(\chiell \big(Q_\ell-\chiellg \cR_\ell\big)(I-Q_0)v_h\big)\big\|_{H^1_k(\Omega_\ell)}
\\ \nonumber
&\hspace{6.3cm}+\big\|(I-\cI_{h}^{\ell})\chiell\cR_\ell(I-Q_0)v_h\big\|_{H^1_k(\Omega_\ell)}\Big]^2\\ \nonumber
&\leq 2\Lambda \sum_{\ell=1}^N\Big[ (1+\mu_\ell)\big(1 + C_{\rm PoU}(k\delta_\ell)^{-1}\big) \big\|\big(Q_\ell-\chiellg\cR_\ell\big) (I-Q_0)v_h\big\|_{H^1_k(\Omega_\ell)}
\\ \label{eq:nomagic1}
&\hspace{8cm}
+\mu_\ell\big\|(I-Q_0)v_h\big\|_{H^1_k(\Omega_\ell)}\Big]^2\\ \nonumber
&\leq 2\Lambda\sum_{\ell=1}^N \Big[ (1+\mu_\ell)\big(1 + C_{\rm PoU}(k\delta_\ell)^{-1}\big) \gamma_\ell  \Ccom
(k\delta_\ell)^{-1} \big(1+ (k\delta_\ell)^{-1}\big) \N{(I-Q_0)v_h}_{L^2(\Omega_\ell)} \\ \nonumber
&\hspace{1.5cm}+ \Big(1+(1+\mu_\ell)\big(1 + C_{\rm PoU}(k\delta_\ell)^{-1}\big) (1+\gamma_\ell\Ccont )\Big)\mu_\ell \N{(I-Q_0)v_h}_{H^1_k(\Omega_\ell)}\Big]^2\\ \nonumber
&\leq 4\Lambda\sum_{\ell=1}^N  \Big[ (1+\mu_\ell)^2\big(1 + C_{\rm PoU}(k\delta_\ell)^{-1}\big)^2 \gamma_\ell^2\Ccom^2
(k\delta_\ell)^{-2} \big(1+ (k\delta_\ell)^{-1}\big)^2 \N{(I-Q_0)v_h}_{L^2(\Omega_\ell)}^2 \\ \nonumber
&\hspace{1.5cm}
+ \Big(1+(1+\mu_\ell)\big(1 + C_{\rm PoU}(k\delta_\ell)^{-1}\big)(1+\gamma_\ell \Ccont)\Big)^2 \mu_\ell^2 \N{(I-Q_0)v_h}_{H^1_k(\Omega_\ell)}^2\Big]\\ \nonumber
&\leq 4\Lambda (1+\mu)^2\big(1 + C_{\rm PoU}(k\delta)^{-1}\big)^2 \gamma^2 \Ccom^2
(k\delta)^{-2} \big(1+ (k\delta)^{-1}\big)^2 \sum_{\ell=1}^N\N{(I-Q_0)v_h}_{L^2(\Omega_\ell)}^2 \\ \nonumber
&\hspace{1cm}
+4\Lambda \Big(1+(1+\mu) \big(1 + C_{\rm PoU}(k\delta)^{-1}\big)(1+\gamma \Ccont)\Big)^2\mu^2 \sum_{\ell=1}^N \N{(I-Q_0)v_h}_{H^1_k(\Omega_\ell)}^2\\ \nonumber
&\leq 4\Lambda^2 (1+\mu)^2 \big(1 + C_{\rm PoU}(k\delta)^{-1}\big)^2\gamma^2 \Ccom^2
(k\delta)^{-2} \big(1+ (k\delta)^{-1}\big)^2 \N{(I-Q_0)v_h}_{L^2(\Omega)}^2 \\ \nonumber
&\hspace{1cm}
+4\Lambda^2 \Big(1+(1+\mu)\big(1 + C_{\rm PoU}(k\delta)^{-1}\big)(1+ \gamma \Ccont )\Big)^2\mu^2  \N{(I-Q_0)v_h}_{H^1_k(\Omega)}^2.
\end{align}
Then, by the two inequalities in \eqref{eq:sigma}, 
\begin{align*}
&\big\|(Q- I)v_h \big\|^2_{H^1_k(\Omega)}\\
&\leq 4\Lambda^2 \Big[ (1+\mu)^2 \big(1 + C_{\rm PoU}(k\delta)^{-1}\big)^2\gamma^2 \Ccom^2
(k\delta)^{-2} \big(1+ (k\delta)^{-1}\big)^2 \sigmat^2\\
&\hspace{2.5cm}+ \Big(1+(1+\mu)\big(1 + C_{\rm PoU}(k\delta)^{-1})\big)(1+\gamma \Ccont)\Big)^2 \mu^2 \sigmao^2\Big]\N{v_h}_{H^1_k(\Omega)}^2.
\end{align*}
The result \eqref{eq:main} then follows by using that $a^2+ b^2 \leq (a+b)^2$ for $a,b\geq 0$. 
\epf

\bre[Why we need $\chiellg$ in the preconditioner for the current proof of Theorem \ref{thm:main1} to work]\label{rem:cutoffs}
Replacing $\sum_{\ell=1}^N
\big(\matrixR_\ell^{\chiell}\big)^T\matrixA_\ell^{-1} \matrixR_\ell^{\chiellg}$ in \eqref{eq:matrixB} by \eqref{eq:group:a} 
corresponds to changing $\chiellg$ to $\chiell$ in the definition of $Q_\ell$ \eqref{eq:Qell}.
The bound \eqref{eq:GSZ} then holds with $\chiellg$ replaced by $\chiell$; i.e., we can bound 
$\big\|(Q_\ell -\chiell \cR_\ell) u_h\big\|_{H^1_k(\Omega_\ell)}$ by $\|u_h\|_{L^2(\Omega_\ell)}$ and a small multiple of $\|u_h\|_{H^1_k(\Omega_\ell)}$ 
(with $u_h=  (I-Q_0)v_h$ in the proof of Theorem \ref{thm:main1}).
However, if $\chiellg$ is not inserted in \eqref{eq:magic1} (via $\chiell= \chiell\chiellg$), then
\beqs
\big\|\big(Q_\ell-\chiellg\cR_\ell\big) (I-Q_0)v_h\big\|_{H^1_k(\Omega_\ell)}
\eeqs 
in \eqref{eq:nomagic1} is replaced by 
\beq\label{eq:nomagic12}
\big\|(Q_\ell -\cR_\ell)(I-Q_0)v_h\big\|_{H^1_k(\Omega_\ell)}.
\eeq
Since $Q_\ell(I-Q_0)v_h$ contains only information about $(I-Q_0)v_h$ on the support of $\chiell$ (via \eqref{eq:Qell}), \eqref{eq:nomagic12} cannot be small. 

Replacing $\sum_{\ell=1}^N
\big(\matrixR_\ell^{\chiell}\big)^T\matrixA_\ell^{-1} \matrixR_\ell^{\chiellg}$ in \eqref{eq:matrixB} by \eqref{eq:group:b} 
corresponds to replacing the definition of $Q_\ell$ \eqref{eq:Qell} by 
\beqs
a_{\ell}\big( Q_\ell u_h, v_{\ell,h} \big)= a\big( u_h, E_{\ell,h} v_{\ell,h}) \big)\quad\tfa v_{\ell,h} \in \subdomain,
\eeqs
where $E_{\ell,h}$ is the ``nodewise extension operator" which extends $v_{\ell,h}$ to an element of $\Hilbert_h$ by setting the values of $E_{\ell ,h} v_{\ell,h}$ at nodes outside $\Omega_\ell$ to be zero (see, e.g., \cite[Equation 2.2]{GoGrSp:23}).
Now, \eqref{eq:nomagic4} in the proof of Lemma \ref{lem:GSZ} becomes 
\beq\label{eq:nomagic6}
a_\ell \big(Q_\ell u_h - \cR_\ell u_h, v_{\ell,h}\big) 
=a\big(u_h, E_{\ell,h} v_{\ell,h})\big) - a_\ell\big(\cR_\ell u_h, v_{\ell,h}).
\eeq
Now, $a$ and $a_\ell$ agree ``in the interior of $\Omega_\ell$" via Assumption \ref{ass:supports}, but they do not necessarily agree on all of $\Omega_\ell$ (since $a_\ell$ might have an absorbing layer near $\partial\Omega_\ell$ or boundary condition on $\partial\Omega_\ell$). Therefore, 
it is not clear how one can obtain a commutator from right-hand side of \eqref{eq:nomagic6} (and hence gain regularity and obtain a bound involving $\|u_h\|_{L^2(\Omega_\ell)}$ as in the proof of Lemma \ref{lem:GSZ}).
\ere

\section{
Assumptions 
\ref{ass:super}, 
\ref{ass:dis},
\ref{ass:supports}, and 
\ref{ass:commutator}
are valid for Helmholtz problems with piecewise-polynomial fine spaces}\label{sec:verify_ass}

\S\ref{sec:ass} contains five main assumptions:~Assumptions \ref{ass:super}, 
\ref{ass:dis},
\ref{ass:supports}, 
\ref{ass:commutator}, and 
\ref{ass:coarse}.
Here we show that the first four of these -- i.e., those not involving the coarse space -- are satisfied for Helmholtz problems with piecewise-polynomial fine spaces.

\subsection{Definition of the finite-element space}\label{sec:space}

We consider a partition of $\Omega$ into a family of conforming meshes $\cT_h$ of
(potentially curved) 
elements $K$. 
For $K \in \cT_h$, let  $h_K$ be the diameter of $K$, and 
let $\LF_K: \widehat K \to K$ be the mapping
between the reference simplex $\widehat K$ and the element $K$. 
Fixing a polynomial degree $p \in\mathbb{Z}^+$, we associate with $\CT_h$ the 
finite element space
\begin{equation*}
\Hilbert_h
:=
\left \{
v_h \in \Hilbert
\; : \;
v_h \circ \MAP_K \in \CP_p(\hK)
\right \}.
\end{equation*}

We consider \emph{both} the case when the element maps are affine, and thus the mesh $\cT^h$ is simplicial \emph{and} the case when the elements are curved.

\begin{assumption}\label{ass:affine}
 $\cT^h$ is  a family of
conforming simplicial meshes of $\Omega$ 
that are 
shape-regular (in the sense of \cite[Chapter 3, Page 111]{Ci:02})
\end{assumption}

\begin{lemma}\label{lem:interpolation_affine}
Given $p\in\mathbb{Z}^+$, if Assumption \ref{ass:affine} holds 
then there exists
$\Cint>0$ and 
 a nodal interpolation operator
$\CI_h: H^2(\Omega) \cap \Hilbert \to \Hilbert_h$ such that, 
for all $\element \in \cT^h$ and 
for all $v \in H^{p+1}(\element ) \cap \Hilbert$,
\begin{equation}\label{eq:interpolation_affine}
h_\element ^{-1} \big\|(I-\CI_h)v\big\|_{L^2(\element )} + \big|(I-\CI_h)v\big|_{H^1(\element )}
\leq
\Cint h_\element ^p |v|_{H^{p+1}(\element )}.
\end{equation}
\end{lemma}

\bpf[References for the proof]
See, e.g., \cite[\S3.1]{Ci:02}, \cite[Theorem 4.4.4]{BrSc:08}.
\epf

\begin{assumption}[Curved elements]
\label{ass:curved}
Let $L$ be the characteristic length scale of $\Omega$. There exists $\Ccurve>0$ such that, for all $h>0$, all $K\in \cT^h$, and $1 \leq |\alpha| \leq p+1$,
\begin{equation}
\label{eq_mapK}
\|\partial^\alpha \MAP_K\|_{L^\infty(\hK)}
\leq
\Ccurve 
L\left (\frac{h_K}{L}\right )^{|\alpha|}
\quad\tand\quad
\|\partial^\alpha(\MAP_K^{-1})\|_{L^\infty(K)}
\leq
\Ccurve h_K^{-|\alpha|}.
\end{equation}
\end{assumption}

Note that the bound \eqref{eq_mapK} with $|\alpha|=1$ implies that the mesh is shape-regular (in the sense of \cite[Chapter 3, Page 111]{Ci:02}).

\bre[Assumption \ref{ass:curved} is satisfied for a piecewise $C^{p+1}$ domain]
Given a piecewise $C^{p+1}$ domain, element maps satisfying Assumption \ref{ass:curved} are constructed in \cite[\S6]{Be:89}, building on the 2-d results of 
\cite{Sc:73, Zl:73} and the isoparametric elements in general dimension of \cite{Le:86}.
\ere

\begin{lemma}[Interpolation on curved meshes \cite{CiRa:72, Le:86, Be:89}]\label{lem:interpolation_curved}
Given $p\in\mathbb{Z}^+$, if Assumption \ref{ass:curved} holds 
then there exists
$\Cint>0$ and 
 a nodal interpolation operator
$\CI_h: H^2(\Omega) \cap \Hilbert \to \Hilbert_h$ such that 
for all $\element \in \cT^h$ and 
for all $v \in H^{p+1}(\element ) \cap \Hilbert$,
\begin{equation}\label{eq:interpolation_curved}
h_\element ^{-1} \big\|(I-\CI_h)v\big\|_{L^2(\element )} + \big|(I-\CI_h)v\big|_{H^1(\element )}
\leq
\Cint h_\element ^p \N{v}_{H^{p+1}(\element )}.
\end{equation}
\end{lemma}

\bpf[References for the proof]
The result under Assumption \ref{ass:curved} is proved in  \cite[Theorem 2]{CiRa:72}; see also \cite[Theorem 1]{Le:86}, \cite[Theorem 5.1]{Be:89}.
\epf

\subsection{Satisfying Assumption \ref{ass:super} (super-approximation)}

\begin{lemma}\label{lem:super}
Let $\Hilbert_h$ be as in \S\ref{sec:space} and 
$\{ \Omega_\ell \}^N_{\ell =1}$ be as in \S\ref{sec:ass2}.
Suppose that \emph{either} Assumption \ref{ass:affine} or Assumption \ref{ass:curved} holds. 
Suppose, additionally, that there exists $C_{\rm PoU}'>0$ such that, for $\ell=1,\ldots,N$, 
 $\chiell, \chiellg \in C^{p,1}(\element )$ with
 \beq\label{eq:PoU_new}
\max\big\{\N{\partial^\alpha \chiell}_{L^\infty(\element )},
\big\|\partial^\alpha \chiellg\big\|_{L^\infty(\element )}\big\}
\leq C_{\rm PoU}' \delta_\ell^{-|\alpha|}
\eeq
for all $\element  \in \cT^h$ all $0\leq |\alpha|\leq p+1$, and some $\{\delta_\ell\}_{\ell=1}^N$ with $\delta_\ell \geq h_\ell := \max_{\element   \subset \overline{\Omega_\ell}}  h_\element $.

Then \eqref{eq:super_approx} holds 
with $\cI_h^\ell=\cI_h$, the interpolant from Lemma \ref{lem:interpolation_affine}/Lemma \ref{lem:interpolation_curved} and there exists $\Csuper>0$ such that
\beq\label{eq:mu_ell}
\mu_\ell:= \Csuper\bigg( \frac{h_\ell}{\delta_\ell} + \bigg(\frac{h_\ell}{\delta_\ell}\bigg)^p \frac{1}{k\delta_\ell}\bigg).
\eeq
\end{lemma}

Similar super-approximation results (albeit not for curved elements) appear in, e.g., \cite[Property A2]{NiSc:74}, \cite[Property A3]{ScWa:77}, \cite[Property A3]{DeGuSc:11}, \cite{Be:99}.

\bpf[Proof of Lemma \ref{lem:super}]
We give the proof under Assumption \ref{ass:curved}; the proof under Assumption \ref{ass:affine} is almost identical (there is no $L^2$ norm on the right-hand side of \eqref{eq:lastday1} below). 
In the proof $C, C', C'',$ and $C'''$ denote constants depending on $p$ and $C_{\rm PoU}'$ whose values may change from equation to equation. 
By \eqref{eq:interpolation_curved},
\beq\label{eq:lastday1}
\big|(I-\CI_h)(\chi_\ell v_h)\big|_{H^1(\element )}
\leq
C h_\element ^p \N{\chi_\ell v_h}_{H^{p+1}(\element )}
\leq C' h_\element^p \Big( | \chiell v_h|_{H^{p+1}(\element)} + \N{\chiell v_h}_{L^2(\element)}\Big).
\eeq
Now, 
\begin{align}\label{eq:clock2}
| \chiell v_h|_{H^{p+1}(\element)} &\leq C \sum_{m=0}^{p+1} |\chiell|_{W^{m,\infty}(\element)} |v_h|_{H^{p+1-m}(\element)}.
\end{align}
Recall that $\cF_K: \hK \to K$. Let $\widetilde{\cF}_K = \cF_k \circ T_{h_K^{-1}}$, where $T_{h_K^{-1}}(x) = h_K^{-1}x$; i.e., $T_{h_K^{-1}}: h_K \hK \to \hK$ via scaling by $h_K^{-1}$. Therefore $\widetilde{\cF}_K: h_K \hK \to K$ and \eqref{eq_mapK} implies that 
\beq\label{eq:clock1}
\|\partial^\alpha \widetilde{\cF}_K\|_{L^\infty(\hK)}
\leq
C 
\quad\tand\quad
\|\partial^\alpha(\widetilde{\cF}_K^{-1})\|_{L^\infty(h_K \hK)}
\leq
C.
\eeq
Since $v_h\circ \cF_K$ is a polynomial of degree $p$, so is $v_h \circ \widetilde{\cF}_K$ and thus $|v_h \circ \widetilde{\cF}_K|_{H^{p+1}(h_K \hK)} =0$. Therefore, by the chain rule and a standard inverse estimate on shape-regular meshes (see, e.g.,  \cite[Theorem 4.76, Page 208]{Sc:98}),
\begin{align*}
\int_K |\partial^\alpha v_h(x)|^2 \, \rd x 
& = \int_K \big|\partial^\alpha \big( v_h \circ \widetilde{\cF}_K \circ\widetilde{\cF}_K^{-1}\big)(x)\big|^2 \, \rd x \\
&\leq C \int_K \sum_{|\beta|\leq |\alpha|} \big|\partial^\beta \big( v_h \circ \widetilde{\cF}_K\big)(\widetilde{\cF}_K^{-1}(x))\big|^2 \, \rd x \\
& \leq C' \sum_{0\leq j \leq \min\{ p, |\alpha|\}} \big| v\circ \widetilde{\cF}_K\big|^2_{H^j (h_K \hK)}\\
&\hspace{-1cm} \leq C'' \bigg( \| v\circ \widetilde{\cF}_K\|^2_{L^2(h_K \hK)} + \sum_{1\leq j\leq \min\{ p, |\alpha|\}} h_K^{-2(j-1)}  \big| v\circ \widetilde{\cF}_K\big|^2_{H^1 (h_K \hK)}\bigg),
\end{align*}
so that 
\begin{align*}
|v_h|_{H^{|\alpha|}(K)}
&\leq C''' \Big( \N{v_h}_{L^2(K)}+h_K^{1- \min \{p,|\alpha|\}} |v_h|_{H^1(K)} \Big).
\end{align*}
Combining this with \eqref{eq:clock2} and \eqref{eq:PoU_new}, we obtain that 
\beq\label{eq:lastday2}
| \chiell v_h|_{H^{p+1}(\element)} \leq C\bigg( \sum_{m=1}^p \delta_\ell^{-m} h_\element^{m-p}\bigg) |v_h|_{H^1(\element)} + C'
\bigg(\sum_{m=0}^{p+1} \delta_\ell^{-m} \bigg)\N{v_h}_{L^2(K)}.
\eeq
The combination of \eqref{eq:lastday1}, \eqref{eq:lastday2}, and the fact that $h_K\leq h_\ell \leq\delta_\ell$ then implies that 
\begin{align}\nonumber
k^{-1}\big|(I-\CI_h)(\chi_\ell v_h)\big|_{H^1(\element )}
&\leq
C \bigg( \frac{h_K}{\delta_\ell} k^{-1} |v_h|_{H^1(K)} + \bigg(\frac{h_K}{\delta_\ell}\bigg)^p \frac{1}{\delta_\ell k} \N{v_h}_{L^2(K)}\bigg)\\
&\leq
C \bigg( \frac{h_K}{\delta_\ell}  + \bigg(\frac{h_K}{\delta_\ell}\bigg)^p \frac{1}{\delta_\ell k}\bigg) \N{v_h}_{H^1_k(K)}.
\label{eq:lastday3}
\end{align}
Arguing in a very similar way, we obtain the inequality
\beq\label{eq:lastday4}
\N{(I-\CI_h)(\chi_\ell v_h)}_{L^2(\element )}
\leq
C\bigg(\frac{h_K}{\delta_\ell}\bigg) \N{v_h}_{L^2(K)},
\eeq
Adding \eqref{eq:lastday3} and \eqref{eq:lastday4}, summing over $K$ intersecting the support of $\chi_\ell$, and using that
$h_K\leq h_\ell$ 
 for such elements, we obtain the bound involving $\chiell$ in \eqref{eq:super_approx} with $\mu_\ell$ given by \eqref{eq:mu_ell}. The proof of the bound involving $\chiellg$ in \eqref{eq:super_approx} is identical.
\epf

\subsection{Satisfying Assumption \ref{ass:supports} ($a$ and $a_\ell$ agree ``in the interior of $\Omega_\ell$ and near $\partial\Omega$")}

\ble\label{lem:supports}
Suppose that $a(\cdot,\cdot)$ is defined by \eqref{eq:sesqui_intro},  $a_\ell(\cdot,\cdot)$ is defined by \eqref{eq:sesquiell_intro},
\beq\label{eq:supports1}
\coeffA_\ell \equiv \coeffA, \quad \coeffB_\ell \equiv \coeffB, \quad \coeffc_\ell \equiv \coeffc \quad \ton \quad\supp\, \cI_h(\chiellg)
\eeq
and
\beq\label{eq:supports2}
\theta\equiv \theta_\ell \quad\ton\quad \supp\, \cI_h(\chiellg) \cap \partial\Omega.
\eeq
Then Assumption \ref{ass:supports} holds.
\ele

\bpf
This follows immediately from the definitions of $a(\cdot,\cdot)$ and $a_\ell(\cdot,\cdot)$. 
\epf

\subsection{Satisfying Assumption \ref{ass:commutator} (the commutator property \eqref{eq:commutator}}

\ble\label{lem:commutator}
Suppose that, for $\ell=1,\ldots,N$,  $a_\ell(\cdot,\cdot)$ is defined by \eqref{eq:sesquiell_intro} with $A_\ell\in C^{0,1}(\Omega_\ell;\Com^{d\times d})$, $B_\ell\in L^\infty(\Omega_\ell;\Com^d)$, $c_\ell^{-2}\in L^\infty(\Omega_\ell,\Com)$ and with their norms bounded independently of $k$, and $\chiellg \in C^{p,1}(\element )\cap C^{1}(\Omega_\ell)$ satisfying 
\eqref{eq:PoU_new}.
If $u(A_\ell\nabla \chiellg) \cdot \nu =0$ on $\partial\Omega$ for all $u\in \cR_\ell(\Hilbert)$,  
then Assumption \ref{ass:commutator} holds.
\end{lemma} 

\bpf
Using the definition of $a_\ell(\cdot,\cdot)$ \eqref{eq:sesquiell_intro}, splitting $\Omega_\ell$ into a sum of mesh elements, and integration by parts, 
we obtain that 
\begin{align*}
&\big| a_\ell(\chiellg u, v)- a_\ell(u,\chiellg v)\big| = \int_{\Omega_\ell} k^{-2} (\coeffA_\ell\nabla \chiellg)\cdot\big(u\overline{\nabla v} - \nabla u \overline{v} \big) + k^{-1}(\coeffB_\ell\cdot \nabla\chiellg)u\overline{v}\\
&\qquad= \int_{\Omega_\ell} k^{-2} (\coeffA_\ell\nabla \chiellg)\cdot u\overline{\nabla v} 
+ k^{-1}(\coeffB_\ell\cdot \nabla\chiellg)u\overline{v}
\\
 &\qquad\quad
+k^{-2} \bigg[
\sum_{K\cap \Omega_\ell \neq \emptyset} \int_{K} u\nabla\cdot\big((\coeffA_\ell\nabla \chiellg)\overline{v}\big) + \int_{\partial K} (A_\ell \nabla\chiellg)\cdot \nu_K \,u\overline{v}\bigg],
\end{align*}
where $\nu_K$ is the outward-pointing unit normal vector to $K$. Since $\chiellg \in C^1(\Omega_\ell)$ the sum over the $\partial K$s that do not touch $\partial \Omega$ is zero. By the assumption that $u(A_\ell\nabla \chiellg) \cdot \nu =0$ on $\partial\Omega$, the sum over the $\partial K$s that touch $\partial\Omega$ is also zero. 
The bound \eqref{eq:commutator} then follows by the bound \eqref{eq:PoU_new}, the definition of $\|\cdot\|_{H^1_k(\Omega_\ell)}$ \eqref{eq:1knorm}, and the assumption that the norms of $A_\ell$ and $B_\ell$ are independent of $k$.
\epf

\bre
The assumption that $\chiellg \in C^1(\Omega_\ell)$ can, in principle, be replaced by $\chiellg \in C^{0,1}(\Omega_\ell)$ together with assumptions on the size of the jumps of the normal derivatives of $\chiellg$ across element boundaries. However, this would substantially complicate the argument of Lemma \ref{lem:commutator}, and we leave it for future work.
\ere

\bre[Constructing $\chiellg$ as in Lemma \ref{lem:commutator}]
The functions $\{\chiellg\}_{\ell=1}^N$ as in Lemma \ref{lem:commutator} -- with, in particular, each $\chiellg \in 
C^{p,1}(\element )\cap C^{1}(\Omega_\ell)$ -- can be constructed by modifying the standard PoU construction appearing in, e.g., \cite[Lemma 5.7]{DoJoNa:15}, \cite[\S3.2]{ToWi:05}. The standard construction defines the PoU functions in terms of the distance function; now one convolves the distance function $p$ times with a piecewise-linear hat function at scale $\delta_\ell$ -- since both the distance function and the hat function are $C^{0,1}$, the resulting function is $C^{p,1}(\Omega_\ell)$. 
\ere

\subsection{Satisfying Assumption \ref{ass:dis} (the discrete inf-sup condition)}

As described in the discussion after Theorem \ref{thm:informal1}, we focus on the situation when $H_\ell\lesssim k^{-1}$ (both for brevity, and because this is the most interesting case).

When verifying Assumption \ref{ass:dis}, 
care is needed if $\partial \Omega$ is disconnected -- i.e., there is an impenetrable obstacle -- and the boundary conditions on the part of $\partial\Omega$ corresponding to the obstacle do not match the boundary conditions on the DD subdomains (e.g., if the obstacle has Neumann boundary conditions, but one wants to use a PML with Dirichlet boundary condition). For brevity, we therefore make the following simplifying assumption (but emphasise that more general situations can be considered).

\begin{assumption}\label{ass:simple}
One of the following holds.

(i) $\Hilbert= H^1_0(\Omega)$ and $\Hilbert_\ell = H^1_0(\Omega_\ell)$.

(ii) $\partial\Omega$ is connected, $\Hilbert= H^1(\Omega)$, and $\Hilbert_\ell = H^1(\Omega_\ell)$.
\end{assumption}

Observe that, in both these cases, $\Hilbert_\ell$ is indeed a closed subspace of $\cR_\ell(\Hilbert)$ -- as required in Section \ref{sec:ass2}.

\bre[Examples of global and local problems satisfying Assumption \ref{ass:simple}]
The choices of $\Hilbert$ and $\Hilbert_\ell$ in Part (i) of Assumption \ref{ass:simple} are used  when the domain $\Omega$ and all the subdomains $\Omega_\ell$ are hyperrectangles, $A$ and $B$ correspond to a Cartesian PML near $\partial\Omega$ and $A_\ell$ and $B_\ell$ correspond to a Cartesian PML near $\partial \Omega_\ell$; these $A, B, A_\ell,$ and $B_\ell$ satisfy \eqref{eq:supports1}, 
with their precise definitions given in, e.g., \cite[\S1.4.1-1.4.2, \S8.1]{GGGLS2}.

 The choices of $\Hilbert$ and $\Hilbert_\ell$ in Part (i) of Assumption \ref{ass:simple} are used when $a(\cdot,\cdot)$ is defined by \eqref{eq:sesqui_intro},  $a_\ell(\cdot,\cdot)$ is defined by \eqref{eq:sesquiell_intro}, with $\theta, \theta_\ell \neq 0$ and satisfying \eqref{eq:supports2}; i.e., impedance boundary conditions are imposed on $\partial \Omega$ and $\partial \Omega_\ell, \ell=1,\ldots,N$, with the global and local impedance parameters (i.e., $\theta$ and $\theta_\ell$, $\ell=1,\ldots,N$) equal to each other in the regions where they are both defined.
\ere

\begin{lemma}[Discrete inf-sup condition for subdomain widths $\lesssim k^{-1}$]\label{lem:dis_verify}
Let $C_1,k_0>0$ and $a_\ell(\cdot,\cdot)$ be given by \eqref{eq:sesquiell_intro}. 
Suppose that there exist $c,s>0$ such that for all $k\geq k_0$ $\Re \coeffA_\ell\geq c>0$ (in the sense of quadratic forms) and $\coeffA_\ell\in C^{1+s}(\Omega_\ell)$, $\Hilbert$ and $\Hilbert_\ell$ satisfy Assumption \ref{ass:simple}, given $F\in (\Hilbert_\ell)^*$, the solution to the variational problem
\beq\label{eq:vpell}
\tfind u_\ell \in \Hilbert_\ell \quad\tst\quad a_\ell(u_\ell,v_\ell)=F(v_\ell)\quad\tfa v_\ell \in \Hilbert_\ell
\eeq
is unique, $\Hilbert_h$ satisfies \emph{either} Assumption \ref{ass:affine} \emph{or} Assumption \ref{ass:curved}, and each subdomain $\Omega_\ell$ is the union of elements of $\cT^h$. 

Then there exists $C_2,C_3>0$ such that if $k\geq k_0$, $\max_\ell (kH_\ell) \leq C_1$, and $\max_\ell (kh_\ell) \leq C_2$, then $\gamma \geq C_3$.
\end{lemma}

\bpf
The bound $\Re \coeffA_\ell\geq c>0$ implies that $a_\ell(\cdot,\cdot)$ satisfies a G\aa rding inequality. The fact that the solution to \eqref{eq:vpell} is unique, combined with Fredholm theory (see, e.g., \cite[Theorem 2.33]{Mc:00}), then implies that the solution to \eqref{eq:vpell} exists.

The change of variables $x= H_\ell \widetilde{x}$ transforms \eqref{eq:vpell} to a Helmholtz problem on an order-one domain with wavenumber $kH_\ell$. 
Since $kH_\ell \leq C_1$, the wavenumber is then 
bounded independently of $k$. The $L^2\to H^1_k$ norm of the solution operator for this problem is therefore bounded independently of $k$. If the solution $u_\ell$ of \eqref{eq:vpell} is in $H^{1+\epsilon}(\Omega_\ell)$, for some $\epsilon>0$, when $F(v_\ell) := \int_{\Omega_\ell} f v_\ell$ for $f\in L^2(\Omega_\ell)$, then the Schatz argument (Appendix \ref{sec:Schatz_app}) combined with the interpolation result \eqref{eq:interpolation_affine}/\eqref{eq:interpolation_curved} and the fact that the subdomains are resolved by the mesh $\cT^h$ imply that the sequence of Galerkin solutions is quasi-optimal, with quasi-optimality constant independent of $k$. The result \cite[Theorem 4.2]{MeSa:10} then implies that $\gamma_\ell \geq C_3$ (using that the solution operator is bounded independently of $k$, and noting that \cite{MeSa:10} work with the weighted norm $\vertiii{\cdot}_{H^1_k}$ defined in the paragraph after \eqref{eq:1knorm}).

It is therefore sufficient to prove that the solution of \eqref{eq:vpell} is in $H^{1+\epsilon}(\Omega_\ell)$, for some $\epsilon>0$, for $L^2(\Omega_\ell)$ data.
If the Dirichlet trace of $u$ is in $H^1(\partial\Omega_\ell)$, then this result with $\epsilon=1/2$ follows from \cite[Theorem 3.1, Equation 3.8]{MiMiTa:01}, using that 
$\coeffA_\ell\in C^{1+s}(\Omega_\ell)$.
If $\Hilbert$ and $\Hilbert_\ell$ satisfy the first condition in Assumption \ref{ass:simple}, then the Dirichlet trace of $u$ on $\partial\Omega_\ell$ is zero, and thus in $H^1(\partial\Omega_\ell)$. 
If $\Hilbert$ and $\Hilbert_\ell$ satisfy the second condition in Assumption \ref{ass:simple}, then $\partial_n u_\ell = \ri k\theta u_\ell$ on $\partial \Omega_\ell$.
If $\theta=0$, $\partial_n u_\ell=0$, and if $\theta\neq 0$, $\partial_n u_\ell\in H^{1/2}(\partial\Omega_\ell)$. In both cases, therefore, $\partial_n u_\ell\in L^2(\partial\Omega_\ell)$. The result 
\cite[\S5.1.2.]{Ne:67} (see also \cite[Theorem 4.24]{Mc:00}) 
then implies that 
the trace of $u$ is in $H^1(\partial\Omega_\ell)$, and the proof is complete.
\epf

\section{Conditions under which Assumption \ref{ass:coarse} (i.e., the assumption on the coarse space) holds}\label{sec:coarse}

\subsection{Sufficient conditions for Assumption \ref{ass:coarse} to hold via the Schatz argument}
Given $f\in L^2(\Omega)$, let $\cS^* f\in\Hilbert$ be the solution of the variational problem 
\beq\label{eq:Csol}
a(v, \cS^* f) =\int_\Omega v \overline{f} \quad\tfa v\in \Hilbert.
\eeq

\ble[Sufficient conditions for Assumption \ref{ass:coarse} to hold via the Schatz argument]\label{lem:coarse}
Let 
\beq\label{eq:adjoint_approx}
\eta(\cV_0):= 
\big\| (I-\Pi_0)\mathcal{S}^*\big\|_{L^2(\Omega)\to H^1_k(\Omega)},
\eeq
where $\Pi_0: \Hilbert\to \coarse$ is the orthogonal projection in the $H^1_k(\Omega)$ norm \eqref{eq:1knorm}. 
Suppose that $Q_0$ is defined by \eqref{eq:Q0}. Then there exists $c>0$ such that if 
\beq\label{eq:sleep3a}
\eta(\cV_0)\leq 
c,
\eeq
then $Q_0: \fine\to \coarse$ is well-defined and such that 
\beq\label{eq:SchatzL2}
\N{(I-Q_0)v_h}_{L^2(\Omega)}\leq \eta(\cV_0) \Ccont \N{(I-Q_0)v_h}_{H^1_k(\Omega)}
\eeq
and 
\beq\label{eq:QO1}
\N{(I-Q_0)v_h}_{H^1_k(\Omega)} \leq 2 \Ccont \big\|(I-\Pi_0) v_h\big\|_{H^1_k(\Omega)}.
\eeq
Therefore (since $\|I-\Pi_0\|_{H^1_k(\Omega)\to H^1_k(\Omega)}\leq 1$) the first inequality in \eqref{eq:sigma} holds with $\sigmat= 2\eta(\cV_0) (\Ccont)^2$ and 
the second inequality in \eqref{eq:sigma} holds with $\sigmao= 2\Ccont$.
\ele

\bpf
This follows from the Schatz argument (recapped as Theorem \ref{thm:Schatz} below). Recall that  the notation $\eta(\cV_0)$ was introduced in \cite[Equation 7]{Sa:06}.
\epf

\subsection{Approximation spaces in the literature satisfying Assumption \ref{ass:coarse} with $\sigma_{L^2}$ small}\label{sec:coarse_other}

\paragraph{Approximation space of \cite{MaAlSc:23}, used as a coarse space in \cite{MaAlSc:24}, following the use of a related coarse space in \cite{HuLi:24} }

The ``discrete MS-GFEM'' of \cite[\S4]{MaAlSc:23} (which creates a multiscale coarse space inside $\Hilbert_h$) satisfies Assumption \ref{ass:coarse}. Indeed, 
\cite[Lemma 3.13/Equation 3.62]{MaAlSc:23} gives conditions under which $\eta$ defined by \eqref{eq:adjoint_approx} is small, with these conditions then carried into the two-level hybrid Schwarz theory in \cite{MaAlSc:24} (see \cite[Theorem 2.17]{MaAlSc:24}). 

The coarse spaces of \cite{HuLi:24} and \cite{MaAlSc:24} are closely related -- see discussion in \cite[\S1]{MaAlSc:24}.
We note that \cite[Equation 6.16]{HuLi:24} is precisely \eqref{eq:SchatzL2} (noting that \cite{HuLi:24} use the $H^1$ norm $\vertiii{\cdot}_{H^1_k(\Omega)}$ defined in the paragraph after \eqref{eq:1knorm}) and the displayed equation immediately after 
 \cite[Equation 6.16]{HuLi:24} is the second inequality in \eqref{eq:sigma}.

\paragraph{Approximation space of \cite{FuLiCrGu:21}, used as a coarse space in \cite{FuGoLiWa:24}}

The wavelet-based multiscale method of \cite{FuLiCrGu:21} satisfies Assumption \ref{ass:coarse} by the bound on $\eta(\Hilbert_0)$ \eqref{eq:adjoint_approx} in \cite[Theorem 4.4]{FuLiCrGu:21}, with conditions for $\eta(\Hilbert_0)$ being small then given in \cite[Equation 4.11]{FuLiCrGu:21}. 
This approximation space is then used as a coarse space in \cite{FuGoLiWa:24}, with the proof of \cite[Proposition 5.1]{FuGoLiWa:24} proving that 
the two inequalities in \eqref{eq:sigma} hold.

\paragraph{Approximation space of \cite{ChHoWa:23}}

The multiscale method of \cite{ChHoWa:23} satisfies Assumption \ref{ass:coarse} by bound on $\eta$ in \cite[Equation 4.4]{ChHoWa:23}.

\paragraph{The localised orthogonal decomposition (LOD) method \cite{Pe:17}, used as a coarse space in \cite{LuXuZhZo:24}}

The multiscale method of \cite{Pe:17} is used as a coarse space in \cite{LuXuZhZo:24}, with \cite[Equation 3.16]{LuXuZhZo:24} and \cite[Equation 3.12]{LuXuZhZo:24} corresponding, respectively, to the first and second bounds in  \eqref{eq:sigma}. 

\subsection{Piecewise polynomials}

Let
\beqs
\Csol:= \sup_{f\in L^2(\Omega)} \frac{ \N{\cS^* f}_{H^1_k(\Omega)}}{\N{f}_{L^2(\Omega)}}.
\eeqs

\ble[Bound on $\eta(\Hilbert_0)$ for piecewise polynomials]\label{lem:eta_pwp}
Suppose that, for $m\in\mathbb{Z}^+$, $\Hilbert_0$ consists of piecewise polynomials of degree $p\leq m$ on a shape-regular mesh with meshwidth $\Hcoarse$ and one of the following holds.

(i) $\Omega$ is $C^{m,1}$, $\coeffA,\coeffB,$ and $\coeffc$ are all $C^{m-1,1}$, and Assumption \ref{ass:curved} holds,

(ii) $\Omega$ is $C^{1,1}$, $\coeffA,\coeffB,$ and $\coeffc$ are all $C^\infty$ and correspond to a radial PML, and Assumption \ref{ass:curved} holds, or 

(iii) $\Omega$ is a convex polygon/polyhedron, $\coeffA,\coeffB,$ and $\coeffc$ are all $C^\infty$ and correspond to a radial complex absorbing potential with, in particular, $\coeffA\equiv I$ and $\coeffB\equiv 0$ near $\partial\Omega$, and Assumption \ref{ass:affine} holds.

Given $k_0>0$ and $M>0$ there exists $C>0$ such that for all $k\geq k_0$,
\beq\label{eq:eta_bound}
\eta(\Hilbert_0)\leq C\big( k\Hcoarse + (k\Hcoarse)^p \Csol + (kL)^{-M}\big).
\eeq
\ele

\bpf
Part (i) (without the term $k^{-M}$ in \eqref{eq:eta_bound}) is proved in \cite{ChNi:20}; see also \cite[Theorem 1.7]{GS3}. 
Part (ii) is proved in \cite[Lemma 2.5]{GLSW1}, and Part (iii) is proved in \cite{IvanEuan1} by adapting the results of \cite[Theorem 1.5]{GLSW1}. 
\epf

\section{Rigorous statement of Corollary \ref{cor:pwp}}\label{sec:6}

The combination of Lemmas \ref{lem:super}, \ref{lem:supports}, \ref{lem:commutator}, and \ref{lem:eta_pwp} give bounds on
$\|\matrixI-\matrixB_L^{-1}\matrixA  \|_{\matrixD_k}$ and $\|\matrixI-\matrixA\matrixB_R^{-1}\|_{\matrixD_k^{-1}}$ when the coarse space consists of piecewise polynomials -- we now summarise these results here. For brevity we only consider 
\bit
\item the first situation in Assumption \ref{ass:simple} ($\Hilbert= H^1_0(\Omega)$ and $\Hilbert_\ell = H^1_0(\Omega_\ell)$),
\item $\theta\equiv \theta_\ell\equiv0$ (i.e., no impedance boundary conditions), and 
\item Case (i) in Lemma \ref{lem:eta_pwp} (covering the most general class of Helmholtz problems under the natural regularity requirements for using degree-$p$ polynomials), 
\eit 
but we emphasise that it is straightforward to write down results about the omitted cases.

\begin{corollary}[Theorem \ref{thm:main1} applied with piecewise-polynomial coarse spaces]\label{cor:pwp2}
\quad\textbf{(The domain and coefficients.)}
For some $m\geq 2$, 
suppose that $\Omega$ is $C^{m,1}$, and $\coeffA,\coeffB,$ and $\coeffc$ are all $C^{m-1,1}$.
Let $\Hilbert= H^1_0(\Omega)$ and let $a(\cdot,\cdot)$ be defined by \eqref{eq:sesqui_intro} with $\theta\equiv 0$. Assume that the solution of the variational problem \eqref{eq:vp} exists and is unique. 
Let $\Csol$ be defined by \eqref{eq:Csol}.

\textbf{(The fine space.)}
Let $\Hilbert_h$ consist of piecewise-polynomials on a shape regular mesh $\cT^h$ with polynomial degree $p\leq m$ with the element maps satisfying Assumption \ref{ass:curved}. 

\textbf{(The subdomains.)}
Let $\{\Omega_\ell\}_{\ell=1}^N$ be Lipschitz domains such that 
\bit
\item each subdomain $\Omega_\ell$ is the union of elements of $\cT^h$,
\item $C_1\leq kH_\ell \leq C_2$ for some $C_1, C_2>0$, 
\item the partition of unity $\{\chiell\}_{\ell=1}^N$ and functions $\{\chiellg\}_{\ell=1}^N$ such that $\chiellg \equiv 1$ on $\supp \chiell$ satisfy  
 $\chiell\in C^{p+1}(\element )\cap C^{0,1}(\Omega_\ell)$, 
  $\chiellg \in C^{p+1}(\element )\cap C^{1}(\Omega_\ell)$  and 
 \beq\label{eq:PoU_new_pwp}
\max\big\{\N{\partial^\alpha \chiell}_{L^\infty(\element )},
\big\|\partial^\alpha \chiellg\big\|_{L^\infty(\element )}\big\}
\leq C_{\rm PoU} k^{|\alpha|}
\eeq
for all $\element  \in \cT^h$ and all $0\leq |\alpha|\leq p+1$.
\eit

\textbf{(The local sesquilinear forms.)}
Let $\Hilbert_\ell = H^1_0(\Omega_\ell)$. Suppose that, for $\ell=1,\ldots,N$,  $a_\ell(\cdot,\cdot)$ is defined by \eqref{eq:sesquiell_intro} with 
\bit
\item $A_\ell\in C^{1+s}(\Omega_\ell;\Com^{d\times d})$ for some $s>0$, $B_\ell\in L^\infty(\Omega_\ell;\Com^d)$, $c_\ell^{-2}\in L^\infty(\Omega_\ell,\Com)$ and with their norms bounded independently of $k$,
\item $\Re \coeffA_\ell\geq c>0$ (in the sense of quadratic forms), and 
\beqs
\coeffA_\ell \equiv \coeffA, \quad \coeffB_\ell \equiv \coeffB, \quad \coeffc_\ell \equiv \coeffc,  \quad\tand \quad\theta\equiv \theta_\ell \quad\ton \supp\, \cI_h(\chiellg).
\eeqs
\eit

\textbf{(The coarse space.)}
$\Hilbert_0\subset \Hilbert_h$ consists of piecewise polynomials of degree $p\leq m$ on a shape-regular mesh with meshwidth $\Hcoarse$, with each coarse-grid element a union of fine-grid elements, and the coarse-space element maps satisfying Assumption \ref{ass:curved}.

\textbf{(The result.)}
Given $m, C_1,C_2, C_{\rm PoU}, s>0$ such that the above hold, and $\epsilon>0$ and $\Ccoarse>1$, there exists $k_0,c>0$ such that if $k\geq k_0$, 
\beq\label{eq:final_threshold}
(k\Hcoarse)^p \Csol \leq c,  \quad\tand\quad h\leq \Hcoarse/\Ccoarse, 
\eeq
then the Galerkin solution \eqref{eq:Galerkin} exists and is quasi-optimal (with quasi-optimality constant independent of $k$) and
\beq\label{eq:final_pwp1}
\max\Big\{\N{\matrixI-\matrixB_L^{-1}\matrixA  }_{\matrixD_k}, \N{\matrixI-\matrixA\matrixB_R^{-1}}_{\matrixD_k^{-1}}\Big\}\leq \epsilon.
\eeq
\end{corollary}

\bpf
The combination of 
Lemmas \ref{lem:super}, \ref{lem:supports}, \ref{lem:commutator}, \ref{lem:dis_verify}, and \ref{lem:coarse}/\ref{lem:eta_pwp} 
(verifying Assumptions \ref{ass:super}, \ref{ass:supports}, \ref{ass:commutator}, \ref{ass:dis}, and \ref{ass:coarse} respectively)
and the fact that $h\leq \Hcoarse/\Ccoarse$ 
imply that \eqref{eq:main} holds with $\gamma$ independent of $k$, $\delta$ proportional to $k^{-1}$, 
\beq\label{eq:reallylastday1}
\mu \leq \frac{h}{\delta}+ C\bigg(\frac{h}{\delta}\bigg)^p \leq C'\big( kh + (kh)^p\big)\leq C''\big(\Csol^{-1/p} + \Csol^{-1}\big),
\eeq
and 
\beq\label{eq:reallylastday2}
\sigmat \leq C \big(k\Hcoarse + (k\Hcoarse)^p \Csol\big).
\eeq
(Note that the requirement in Lemma \ref{lem:commutator} that $\chiellg$ is such that $u(A_\ell\nabla \chiellg) \cdot \nu =0$ on $\partial\Omega$ for all $u\in \cR_\ell(\Hilbert)$ is satisfied via the zero Dirichlet boundary condition on $\partial\Omega$.)

Recall that $\Csol \geq C kL$. 
The bound 
\beq\label{eq:Wed1}
\|(I-Q)v_h\|_{H^1_k(\Omega)}\leq \epsilon \N{v_h}_{H^1_k(\Omega)} \tfa v_h\in \Hilbert_h
\eeq
 follows from \eqref{eq:main} since $\sigmat$ is made sufficiently small by \eqref{eq:final_threshold} (compare \eqref{eq:final_threshold} and \eqref{eq:reallylastday2}) and $\mu \to 0$ as $k\to \infty$ (via \eqref{eq:reallylastday1}). The bound \eqref{eq:final_pwp1} then follows from \eqref{eq:Wed1}  by 
Appendix \ref{sec:appendix} and \eqref{eq:ip2}.

The result about quasi-optimality of the Galerkin solution follows from the Schatz argument (Appendix \ref{sec:Schatz_app}) and the analogue of Lemma \ref{lem:eta_pwp} with the coarse space replaced by the fine space. 

We obtain the result about right preconditioning, via \eqref{eq:adjoint2a}, by showing that 
the assumptions of Theorem \ref{thm:main1} are satisfied for the adjoint sesquilinear form.
Lemmas \ref{lem:super}, \ref{lem:supports}, \ref{lem:commutator}, \ref{lem:dis_verify} 
(verifying Assumptions \ref{ass:super}, \ref{ass:supports}, \ref{ass:commutator}, and \ref{ass:dis}, respectively)
hold immediately for the adjoint sesquilinear form. 
To apply Lemma \ref{lem:eta_pwp} (for Assumption \ref{ass:coarse}) we observe that  $\|\cS^*\|_{L^2\to L^2}$ = $\|\cS\|_{L^2\to L^2}$,
and thus the $L^2\to H^1_k$ norms of these operators have the same dependence on $kL$ (since an $L^2 \to L^2$ bound implies an $L^2 \to H^1_k$ bound by Green's identity; see, e.g., \cite[Lemmas 3.10 and A.10]{GrPeSp:19}). The combination of Lemmas \ref{lem:eta_pwp} and \ref{lem:coarse} therefore implies that the condition \eqref{eq:final_threshold} ensures that Assumption \ref{ass:coarse} holds for the adjoint problem, and the proof is complete.
\epf

\bre[The numerical experiments in \cite{BoDoJoTo:21} for piecewise-polynomial coarse spaces]
The experiments in \cite{BoDoJoTo:21} consider discretisations with degree-two polynomials, a fixed number of points per wavelength, and an additive Schwarz preconditioner rather than a hybrid preconditioner. Nevertheless, the results of these experiments are somewhat consistent with Corollary \ref{cor:pwp2} in that they show that the number of GMRES iterations 
\bit
\item grows slowly with $k$ when the subdomain widths $\sim k^{-1}$, 
\item grows with $k$ if the coarse space does not resolve the oscillatory/propagative nature of the solution.
\eit
In more detail:~the \emph{grid coarse space method} of \cite{BoDoJoTo:21} involves FEM discretisations with  fine and coarse polynomial degrees equal to two, 
10 points per wavelength in the fine space, and 5 points per wavelength in the coarse space (i.e., both $h$ and $H_{\rm coarse} \sim k^{-1}$)
 and GMRES is then applied with an additive Schwarz preconditioner with impedance boundary conditions on the subdomains and minimal overlap -- we expect the hybrid preconditioner with generous overlap to have fewer GMRES iterations than in this set up.
When $k$ is doubled and the number of subdomains increases by $2^d$ (so that the number of degrees of freedom per subdomain is kept constant, the number of iterations 
goes from $41\, (f=10, N=40)$ to $44\, (f=20, N=160)$ in \cite[Table 1]{BoDoJoTo:21} for the 2-d Marmousi model
and from $11\, (k=100, N=20)$ to $16\, (k=200, N=160)$ \cite[Table 7]{BoDoJoTo:21} for the 3-d cobra cavity.
 Furthermore, \cite[Table 9]{BoDoJoTo:21} shows that the number of iterations is large if 
 there are only 5 points per wavelength in the fine space, and 2.5 points per wavelength in the coarse space.
 \ere

\begin{appendix}

\section{The matrix form of the operator $Q$ \eqref{eq:Q}}\label{sec:appendix}

\paragraph{Additional notation for the fine space $\fine$}

Denote the nodes of $\cT^h$ by 
$\cN^h = \{x_j:j \in \cJ_h\}$, where $\cJ_h$ is  a suitable index
set. 
Let 
$\{ \phi_j : j \in \cJ_h \}$ 
be the standard nodal basis for $\cV^h$.
Let 
\beq\label{eq:matrixA}
(\matrixA)_{ij}:= a( \phi_j, \phi_i) \quad \tfor i,j \in \cJ_h,
\eeq
so that the Galerkin equations 
\beq\label{eq:Galerkin}
\tfind u_h\in \Hilbert_h\tst a(u_h,v_h) =F(v_h)\quad\tfa v_h\in \Hilbert_h.
\eeq
are equivalent to the linear system 
$\matrixA \bu =\bff.$

\paragraph{Restriction matrices on the fine grid}

Denote the freedoms for $\cV_{\ell,h}$ by 
$\cN^h(\Omega_\ell) = \{x_j:j \in \cJ_h(\Omega_\ell)\}$, where $\cJ_h(\Omega_\ell)$ is a suitable
index set. The  
basis for
$\cV_{\ell,h}$ 
can then be  written as $\{ \phi_j : j \in \cJ_h(\Omega_\ell)
\}$. Let 
\beqs
(\matrixA_\ell)_{ij}:= a_\ell( \phi_j, \phi_i)\quad\tfor i,j \in \cJ_h(\Omega_\ell);
\eeqs
i.e., $\matrixA_\ell$ is the Galerkin matrix of $a_\ell(\cdot,\cdot)$. 
%
%
Let 
\beq\label{eq:weightedR}
(\matrixR_\ell^{\chiell})_{jj'} := \delta_{jj'}\chiell(x_j) \,\,\tand\,\, (\matrixR_\ell^{\chiellg})_{jj'} := \delta_{jj'}\chiellg(x_j), \quad j \in \cJ_h(\Omega_\ell), j'\in \cJ_h.
\eeq

\paragraph{The coarse grid and associated restriction matrices}

Let  $\{\cT^{H}\}$ be  a sequence of shape-regular, simplicial meshes on 
$\overline{\Omega}$, with mesh diameter $H$. We assume that each element of $\cT^H$ consists of the union of a set of 
fine-grid elements. 
Let $\cJ_H$ be the set of 
 coarse mesh nodes, so that 
 $\{\Phi_p : p\in \cJ_H\}$ is the nodal basis for $\cV_0$. 
Since $\coarse\subset \fine$, there exists a matrix $\matrixR_0$ such that 
\beq\label{eq:Phi1}
\Phi_p = \sum_j (\matrixR_0)_{pj} \phi_j. 
\eeq
Let 
\begin{equation}\label{eq:coarsegrid}
\matrixAzero := \matrixR_0 \matrixA \matrixR_0^T.
\end{equation}
In fact, since each element of $\cT^H$ consists of the union of a set of 
fine-grid elements, if $\cI_h$ is the nodal interpolation operator then
\beq\label{eq:Phi2}
\Phi_p =\cI_h \Phi_p= \sum_j \Phi_p(x_j^h) \phi_j,
\eeq
and thus
\begin{equation}\label{eq:restriction}
(\matrixR_0)_{pj}  := \Phi_p(x_j^h) ,  \quad j \in \cJ_h , \quad
p \in \cJ_{H}.
\end{equation}

\ble
$(\matrixAzero)_{pq}= a(\Phi_q,\Phi_p)$; i.e., $\matrixAzero$ is the Galerkin matrix for the variational problem \eqref{eq:Galerkin} 
discretised in $\cV_0$ using the basis $\{\Phi_p: p \in
\cJ_{H}\}$.  
\ele 

\bpf
By the definition \eqref{eq:Phi1} of $\matrixR_0$ and the definition \eqref{eq:matrixA} of $\matrixA$, 
\beqs
a(\Phi_q,\Phi_p)=\sum_j \sum_i (\matrixR_0)_{qj} a(\phi_j, \phi_i) (\matrixR_0)_{pi}= \sum_j \sum_i (\matrixR_0)_{pi} (\matrixA)_{ij} (\matrixR_0^T)_{jq}= (\matrixR_0 \matrixA \matrixR_0^T)_{pq},
\eeqs
and the result follows from \eqref{eq:coarsegrid}.
\epf

\paragraph{The matrix form of the $H^1_k(\Omega)$ inner product}

Let 
\begin{equation}
(\matrixS)_{\ell,m} = \int_{\Omega} \nabla \phi_\ell \cdot \nabla
\phi_m,\qquad
(\matrixM)_{\ell,m} = \int_{\Omega}
\phi_\ell  \phi_m , 
\quad\tand\quad
 \matrixD_k := k^{-2}\matrixS + \matrixM.
\label{eq:weight}
 \end{equation}
 It then follows that 
if $v_h, w_h \in\fine$ with coefficient vectors $\bV, \bW$ then 
\begin{equation}\label{eq:ip}
(v_h, w_h )_{H^1_k(\Omega)}\  =\ \langle \bV, \bW\rangle_{\matrixD_k} . 
\eeq

\paragraph{The matrix form of the operators $Q_{\ell}$}

The fact that the matrix form of $Q$ is $\matrixB_L^{-1}\matrixA$, with $\matrixB_L^{-1}$ defined by \eqref{eq:matrixB}, 
is an immediate consequence of the following result combined with the definition of $Q$ \eqref{eq:Q}.

\begin{theorem}
\label{thm:repQ}
Let $v_h = \sum_{j\in \cJ_h} V_j \phi_j \ \in \fine$. Then, for $\ell=1,\ldots,N$,
\beqs
\cI_h(\chiell Q_{\ell} v_h)  = \sum_{j \in \cJ_h}
\left(\big(\matrixR_\ell^{\chiell}\big)^T\matrixA_\ell^{-1} \matrixR_\ell^{\chiellg} \matrixA \mathbf{V}\right)_j \phi_j,\quad
%
Q_{0} v_h  = \sum_{j \in \cJ_h} \left(\matrixR_0^T\matrixAzero^{-1} \matrixR_0 \matrixA\mathbf{V}\right)_j \phi_j .
\eeqs
\end{theorem}

\bpf
The proof of the first expression is very similar to the proof of  \cite[Theorem 2.10]{GSZ1} (where only one type of weighted restriction matrix is used, in contrast to the two used here). 
The second expression is proved in \cite[Theorem 5.4]{GSV1}. 
\epf

\section{Recap of the Schatz argument and Aubin--Nitsche lemma}\label{sec:Schatz_app}

\vspace{-0.5cm}

\begin{theorem}[The Schatz argument \cite{Sc:74, ScWa:96}]\label{thm:Schatz}
Suppose that the sequilinear form $b: \cH\times\cH\to \Com$ is continuous, i.e., 
\beq\label{eq:contS}
\big|b(u,v)\big| \leq C_{\rm cont} \N{u}_{\cH} \N{v}_{\cH} \quad\tfa u, v\in \cH
\eeq
and satisfies the G\aa rding equality 
\beq\label{eq:GS}
\Re b(v,v)\geq C_{G1}\N{v}_{\cH}^2 -  C_{G2}\N{v}_{\cH_0}^2 \quad\tfa v\in \cH
\eeq
where $\cH_0\subset \cH$. 
Suppose that $\cH_h$ is a finite-dimensional subspace of $\cH$. 
Given $u\in\cH$ and $u_h\in \cH_h$ such that 
\beq\label{eq:GogS}
b(u-u_h,v_h)=0 \quad\tfa v_h\in \cH_h,
\eeq
\beq\label{eq:Schatzepsilon}
\text{ if } \,\,\N{u-u_h}_{\cH_0}\leq \frac{C_{G1}}{\sqrt{2 C_{G2}}} \N{u-u_h}_{\cH}
\,\,\text{ then } \,\,
\N{u-u_h}_{\cH} \leq 2C_{\rm cont} \min_{v_h\in \cH_h} \N{u-v_h}_{\cH}.
\eeq
\end{theorem}

\bpf
By \eqref{eq:GS}, \eqref{eq:GogS}, and \eqref{eq:contS}, for all $v_{h}\in \cH_h$,
\begin{align*}
C_{G1}\N{u-u_h}^2_{\cH} &\leq \Re b\big( u-u_h,u-u_h\big) + C_{G2}\N{u-u_h}^2_{\cH_0}\\
&\leq \Re b\big( u-u_h,u-v_h\big) + C_{G2}\N{u-u_h}^2_{\cH_0}\\
&\leq C_{\rm cont}\N{u-u_h}_{\cH} \N{u-v_h}_{\cH} + C_{G2} \N{u-u_h}^2_{\cH_0},
\end{align*}
and the result  \eqref{eq:Schatzepsilon} follows.
\epf

\ble\emph{(The Aubin--Nitsche lemma \cite[Theorem 3.2.4]{Ci:02})}
\label{lem:eta}
Under the assumptions of Theorem \ref{thm:Schatz}, given $f\in \cH_0$, let $\mathcal{S}^* f$ be the solution of the variational problem
\beq\label{eq:adjoint}
b(w, \cS^* f) = (w,f)_{\cH_0} \quad\tfa w\in \cH.
\eeq
Let 
\beq\label{eq:eta}
\eta(\cH_h):= \sup_{f\in \cH_0} \min_{v_{h}\in \cH_h} \frac{
\N{\mathcal{S}^* f - v_{h}}_{\cH}
}{
\N{f}_{\cH_0}
}.
\eeq
Then 
\beq\label{eq:QOS}
\N{u-u_h}_{\cH_0} \leq \Ccont \eta(\cH_h) \N{u-u_h}_{\cH}.
\eeq
\ele

\bpf
By \eqref{eq:adjoint}, \eqref{eq:GogS}, \eqref{eq:contS}, and \eqref{eq:eta}, for all $v_h\in \cH_h$,
\begin{align*}
\N{u-u_h}^2_{\cH_0}= b\big( u-u_h, \mathcal{S}^* (u-u_h)\big) &= b\big( u-u_h, \mathcal{S}^* (u-u_h)- v_{h}\big) \\
&\leq C_{\rm cont} \N{u-u_h}_{\cH}\N{\mathcal{S}^* (u-u_h) - v_{h}}_{\cH},
\end{align*}
and the result \eqref{eq:QOS} follows.
\epf

\end{appendix}

\section*{Acknowledgements}

The authors thank Théophile Chaumont-Frelet (INRIA, Lille) and Pierre-Henri Tournier (CNRS, Paris) for useful discussions, and the three anonymous referees for their careful reading of the paper and their numerous constructive comments. JG was 
supported by EPSRC grants EP/V001760/1 and EP/V051636/1 and Leverhulme Research Project Grant RPG-2023-325. JG and ES were both supported by the ERC synergy grant ``PSINumScat" 101167139.

\bibliographystyle{siam}
\bibliography{combined.bib}

\end{document}